\documentclass[12pt]{amsart}
\numberwithin{equation}{section}

\begin{document}
\vsize=9in
\oddsidemargin=0cm
\evensidemargin=0cm


\title[Equivalences of
real submanifolds in complex space]{Equivalences of
real submanifolds in complex space}
\author[M. S. Baouendi, L. P. Rothschild, D.
Zaitsev]{M. S. Baouendi, Linda Preiss Rothschild and
Dmitri Zaitsev}
\address{Department of Mathematics, 0112, University
of California at San Diego, La Jolla, CA 92093-0112, USA}
\address{Mathematisches Institut, Universit\"at T\"ubingen,
Auf der Morgenstelle 10, 72076 T\"ubingen, GERMANY}
\email{sbaouendi@ucsd.edu, lrothschild@ucsd.edu,
dmitri.zaitsev@uni-tuebingen.de}
\begin{thanks}{2000 {\em Mathematics Subject
Classification.} 32H02, 32V40, 32V35.}\end{thanks}
\begin{thanks}{The first and second authors are
partially supported by National Science Foundation
grant DMS 98-01258 .}\end{thanks}
\begin{abstract} We show that for any real-analytic
submanifold $M$ in $\Bbb C^N$ there is a proper real-analytic
subvariety $V$ contained in $M$ such that for any $p \in
M\setminus V$, any real-analytic submanifold $M'$ in $\Bbb
C^N$, and any $p' \in M'$, the germs of the submanifolds
$M$ and $M'$ at $p$ and $p'$ respectively are formally
equivalent if and only if they are biholomorphically
equivalent.  More general results for $k$-equivalences are
also stated and proved.
\end{abstract}

\def\Label#1{\label{#1}}
\def\1#1{\ov{#1}}
\def\2#1{\widetilde{#1}}
\def\3#1{\mathcal{#1}}
\def\4#1{\widehat{#1}}

\def\s{s}
\def\k{\kappa}
\def\ov{\overline}
\def\span{\text{\rm span}}
\def\tr{\text{\rm tr}}
\def\xo {{x_0}}
\def\Rk{\text{\rm Rk}}
\def\sg{\sigma}
\def \emxy{E_{(M,M')}(X,Y)}
\def \semxy{\scrE_{(M,M')}(X,Y)}
\def \jkxy {J^k(X,Y)}
\def \gkxy {G^k(X,Y)}
\def \exy {E(X,Y)}
\def \sexy{\scrE(X,Y)}
\def \hn {holomorphically nondegenerate}
\def\hyp{hypersurface}
\def\prt#1{{\partial \over\partial #1}}
\def\det{{\text{\rm det}}}
\def\wob{{w\over B(z)}}
\def\co{\chi_1}
\def\po{p_0}
\def\fb {\bar f}
\def\gb {\bar g}
\def\Fb {\ov F}
\def\Gb {\ov G}
\def\Hb {\ov H}
\def\zb {\bar z}
\def\wb {\bar w}
\def \qb {\bar Q}
\def \t {\tau}
\def\z{\chi}
\def\w{\tau}
\def\Z{\zeta}

\def \T {\theta}
\def \Th {\Theta}
\def \L {\Lambda}
\def\b {\beta}
\def\a {\alpha}
\def\o {\omega}
\def\l {\lambda}

\def \im{\text{\rm Im }}
\def \re{\text{\rm Re }}
\def \Char{\text{\rm Char }}
\def \supp{\text{\rm supp }}
\def \codim{\text{\rm codim }}
\def \Ht{\text{\rm ht }}
\def \Dt{\text{\rm dt }}
\def \hO{\widehat{\mathcal O}}
\def \cl{\text{\rm cl }}
\def \bR{\mathbb R}
\def \bC{\mathbb C}
\def \C{\mathbb C}
\def \bL{\mathbb L}
\def \bZ{\mathbb Z}
\def \bN{\mathbb N}
\def \scrF{\mathcal F}
\def \scrK{\mathcal K}
\def \scrM{\mathcal M}
\def \cR{\mathcal R}
\def \scrJ{\mathcal J}
\def \scrA{\mathcal A}
\def \scrO{\mathcal O}
\def \scrV{\mathcal V}
\def \scrL{\mathcal L}
\def \scrE{\mathcal E}
\def \hol{\text{\rm hol}}
\def \aut{\text{\rm aut}}
\def \Aut{\text{\rm Aut}}
\def \J{\text{\rm Jac}}
\def\jet#1#2{J^{#1}_{#2}}
\def\gp#1{G^{#1}}
\def\gpo{\gp {2k_0}_0}
\def\emmp {\scrF(M,p;M',p')}
\def\rk{\text{\rm rk}}
\def\Orb{\text{\rm Orb\,}}
\def\Exp{\text{\rm Exp\,}}
\def\Span{\text{\rm span\,}}
\def\d{\partial}
\def\D{\3J}
\def\pr{{\rm pr}}
\def\Re{{\rm Re}\,}
\def\Im{{\rm Im}\,}

\newtheorem{Thm}{Theorem}[section]
\newtheorem{Def}[Thm]{Definition}
\newtheorem{Cor}[Thm]{Corollary}
\newtheorem{Pro}[Thm]{Proposition}
\newtheorem{Lem}[Thm]{Lemma}
\newtheorem{Rem}[Thm]{Remark}
\newtheorem{Thmp}[Thm]{Theorem \cite{artin1}}
\maketitle

\tableofcontents
\newpage
\section{Introduction}

     A {\it formal map}
$H: (\bC^N,p)\to (\bC^N,p')$, with $p$ and $p'$ in
$\bC^N$, is a
$\bC^N$-valued formal power series
$$H(Z)=p' + \sum_{|\alpha| \ge 1} a_\alpha
(Z-p)^\alpha,\ \ \
\ a_\alpha \in \bC^N, \quad Z=(Z_1,\ldots,Z_N).
$$
The map $H$ is {\it
invertible} if there exists a formal map
$H^{-1}: (\bC^N, p')\to (\bC^N,p)$ such that
$H(H^{-1}(Z)) \equiv H^{-1}(H(Z))\equiv Z$ (which is equivalent
to the nonvanishing of the Jacobian of $H$ at
$p$).  Suppose $M$ and $M'$ are
real-analytic submanifolds in $\bC^N$ of the same
dimension  given by real-analytic (vector valued)
local defining functions
$\rho(Z,
\bar Z)$ and
$\rho'(Z, \bar Z)$ near  $p\in M$ and $p'\in M'$
respectively. A formal invertible map $H$ as above is
called a {\it formal equivalence} between (the germs) $(M,p)$ and
$(M',p')$ if
$$\rho'\big(H(Z(x)),
\overline {H(Z(x))}\big) \equiv 0$$
in the sense of formal power series in $x$
for some (and
hence for any) real analytic parametrization
$x\mapsto Z(x)$ of $M$ near $p=Z(0)$. If, in
addition,
$H$ is convergent, we say that $H$ is a {\it
biholomorphic equivalence} between $(M,p)$ and
$(M',p')$.  More generally, for any integer
$k>1$, we call a formal invertible mapping
$H:(\bC^N,p)\to (\bC^N,p')$ a $k$-{\em equivalence}
between
$(M,p)$  and $(M',p')$ (see Lemma~\ref{equivss} for equivalent definitions) if
$$\rho'\big(H(Z(x)),
\overline {H(Z(x))}\big) = O(|x|^k).$$
Hence a formal
invertible map $H$ is a formal equivalence
between
$(M,p)$ and $(M',p')$ if and only if it is a
$k$-equivalence for every $k>1$.

If $M$ and $M'$ are as above, we shall say that
$(M,p)$ and
$(M',p')$ are {\it formally equivalent} (resp.\ {\it
biholomorphically equivalent} or
$k$-{\it equivalent}) if there exists a formal
equivalence (resp.\ biholomorphic  equivalence or
$k$-equivalence)  between $(M,p)$ and $(M',p')$.
It is known (see below) that there exist pairs
$(M,p)$ and $(M',p')$, with $M$ and
$M'$ real-analytic submanifolds of $\C^N$,
which are formally equivalent but not
biholomorphically equivalent.
However, as our main result shows, for ``most'' points
$p\in M$, the notions of formal and biholomorphic
equivalences coincide.
More precisely, we prove the following.
\begin{Thm}\Label{2}
        Let $M\subset \bC^N$ be a connected real-analytic
submanifold.  Then there exists a closed proper
real-analytic subvariety $V\subset M$ such that for
every $p\in M\setminus V$,
every real-analytic submanifold $M'\subset
\bC^N$, every $p'\in
M'$, and every integer
$\k>1$, there exists an integer
$k>1$ such that if
$H$ is a
$k$-equivalence between $(M,p)$ and $(M',p')$ then
there exists a biholomorphic equivalence $\4H$
between $(M,p)$ and $(M',p')$ with
$\4H(Z)=H(Z) + O(|Z-p|^\k)$.
\end{Thm}

In fact, a real-analytic subvariety $V\subset M$, for
which Theorem~\ref{2} holds
will be explicitly described in \S\ref{local} below.
An immediate consequence of Theorem~\ref{2} is
the following corollary.

\begin{Cor}\Label{main}
Let $M\subset \bC^N$ be a connected
real-analytic submanifold, $V\subset M$ the
real-analytic subvariety given by Theorem~{\rm\ref{2}}, and
$p\in M\setminus V$.
Then for every real-analytic submanifold $M'\subset
\bC^N$, and every $p'\in
M'$, the following are equivalent:
\begin{enumerate}
\item [(i)] $(M,p)$ and $(M',p')$ are $k$-equivalent
for all
$k>1$;
\item [(ii)] $(M,p)$ and $(M',p')$ are formally
equivalent;
\item [(iii)] $(M,p)$ and $(M',p')$ are
biholomorphically equivalent.
\end{enumerate}
\end{Cor}

It should be noted that in general, the integer $k$ in Theorem~\ref{2} must be
chosen bigger than
$\k$. For example, if $M=M'=\{Z=(z,w)\colon\Im w = |z|^2\}\subset\C^2$, 
one can easily check that
the map $H(z,w)\colon=(z,w+w^{3})$ is a $4$-equivalence between $(M,0)$ and
$(M',0)$. However,  there is no biholomorphic
equivalence $\4H$ between
$(M,0)$ and $(M',0)$  such that $\4H(Z)-H(Z) = O(|Z|^4)$. Indeed, it is known that any
biholomorphic equivalence $\4H$ between
$(M,0)$ and $(M',0)$ 
that differs from the identity (and hence from $H$) by 
$O(|Z|^3)$ must be the identity (see
\cite{CM}), and hence necessarily $\4H(z,w) - H(z,w) = -w^3$. 
   This proves that for this example if
$\k = 4$, one cannot take $k =4$.

The problem of formal versus biholomorphic
equivalence has been studied by a number of
mathematicians. It has been known since the
fundamental work of Chern and Moser
\cite {CM} that if $M$ and $M'$ are real analytic
      hypersurfaces in $\bC^N$
which are Levi nondegenerate at $p$ and $p'$
respectively, then $(M,p)$ and $(M',p')$ are formally
equivalent if and only if they are biholomorphically
equivalent. It should be mentioned here that
Theorem~\ref{2} and its corollary are new
even in the case of a hypersurface. (See also
Remark \ref{hyp} below.)  Although it had been
known (e.g.\ in dynamical systems) that there exist
pairs of structures  which are formally equivalent
(in an appropriate sense) but not biholomorphically
equivalent, to our knowledge the first examples
of  pairs $(M,p)$ and $(M',p')$ of germs of real analytic
submanifolds in
$\bC^N$ which are formally equivalent but not
biholomorphically equivalent are due to Moser and
Webster \cite {MW}. The examples in that paper
consist of real-analytic surfaces $M$ and
$M'$ in
$\bC^2$ with isolated ``complex tangent" at $p$ and
$p'$ respectively. (It is fairly easy to prove
Theorem \ref {2} above in the case of real-analytic
surfaces in $\bC^2$, since outside a real-analytic set
such a surface is either totally real or complex.) The
work of \cite {MW} also contains positive
results for surfaces in $\bC^2$, i.e.\ cases in which
formal and biholomorphic equivalence coincide at
some complex tangent points. We should also mention
further related work by X. Gong \cite{Gong} as well as
recent work by Beloshapka \cite{Belo} and Coffman \cite
{Cof}.

In other recent work of the first two authors jointly
with Ebenfelt \cite{Asian}, \cite{BER99} and
\cite{Conv},  it has been shown that there are
many classes of pairs
$(M,p)$ and
$(M',p')$, where
$M$ and $M'$ are real-analytic generic submanifolds of
$\bC^N$, for which any formal equivalence is
necessarily convergent (see also Corollary~\ref{converge}).  
In particular it follows that
the notions of formal equivalence and
biholomorphic equivalence for such pairs
coincide.
The present paper
treats the more general case where nonconvergent formal equivalences may exist
between $(M,p)$ and $(M',p')$.  Given such a formal equivalence $H$, Theorem~\ref{2} implies the existence of 
a possibly different biholomorphic equivalence 
that coincides with $H$ up to an arbitrarily high preassigned order.
For instance, any formal power series in one variable
of the form $\sum_{j=1}^\infty a_jz^j$, $a_1\not= 0$,
$a_j \in \bR$, may be regarded as a formal equivalence
between $(\bR,0)$ (considered as a germ of a
real submanifold in $\bC$) and itself . By truncating this power series
to any order, one obtains a biholomorphic equivalence
which agrees with the formal equivalence to that order.
 
The organization of the paper is as follows.  In \S \ref
{local} through \S \ref{red} the variety $V$ is
constructed, and a local description of $M$ near a point
$p \in M\setminus V$ is given. The proof of Theorem
\ref{2} is then reduced to the case where $M$ and $M'$
are generic submanifolds which are finitely nondegenerate
at $p$ and $p'$ respectively.  In \S \ref{rings}
through
\S \ref{final} we prove Theorem \ref{2} in that
case. For the proof, we first obtain a universal
parametrization of $k$-equivalences between
$(M,p)$ and $(M',p')$. The construction of this
parametrization is in the spirit of the one in
\cite {BER99} for formal equivalences;
however, the approach used here is somewhat
different and deals with  more general situations.
The main difference is due to the fact that the parametrization
is obtained in terms of finite order jets {\em along a certain 
submanifold} rather than in terms of {\em single} jets at a point $p$.
    From this parametrization we obtain a system of
real analytic equations in the product of 
the above submanifold and the space of jets,
whose exact solutions correspond to biholomorphic equivalences
and whose approximate solutions of finite order correspond to $k$-equivalences.
For this system we apply  approximation theorems
due to Artin \cite {artin1}, \cite {artin2} and a
variant due to Wavrik \cite {Wav}. The proof is
then completed in \S \ref{final}. We
conclude the paper in \S \ref{rfinal} by giving a version of
Corollary \ref{main} for CR maps between CR submanifolds.

The authors wish to thank Leonard Lipshitz for
pointing out to us the article \cite {Wav}, as
well as his related joint work with Denef
\cite{DL}.

\section{Construction of the real
subvariety $V$}\Label{local}

       {}For the remainder of this paper $M$ and $M'$ will
always denote connected real-analytic submanifolds
of $\C^N$ of the same dimension. For any
$p\in M$, we shall define three nonnegative integers:
$r_1(p)$, the  excess codimension of
$M$ at $p$,
$r_2(p)$, the degeneracy of $M$ at $p$,
and $r_3(p)$, the orbit codimension of $M$ at $p$.
We shall show that  these integers  reach their
minima outside  proper real-analytic subvarieties
$V_1,V_2,V_3\subset M$ respectively and shall prove
Theorem~\ref{2} for $V:=V_1\cup V_2\cup V_3$.

Let $M$ be as above, $d$ be the codimension of $M$ in
$\C^N$, and $p_0\in M$ be fixed. Recall
that a (vector valued) local defining function
$\rho=(\rho^1,\ldots,\rho^d)$ near
$p_0$ is  a collection of real valued real-analytic
functions defined in a neighborhood of $p_0$ in
$\C^N$ such that $M=\{Z:\rho(Z,\1Z)=0\}$  near
$p_0$ and
$d\rho^1\wedge\ldots\wedge d\rho^d\ne 0$.
We associate to $M$ a complex submanifold
$\3M\subset\C^{2N}$ (called the {\em
complexification of $M$}) in a neighborhood of
$(p_0,\1p_0)$ in $\bC^N\times \bC^N$ defined by
$\3M:=\{(Z,\Z):\rho(Z,\Z)=0\}$.
Observe that a point $p\in\C^N$ is in $M$
if and only if $(p,\1p)\in \3M$.
We also note that, if
$\2\rho=(\2\rho^1,\ldots,\2\rho^d)$
is another local defining function of $M$ near $p_0$,
then $\2\rho(Z,\1Z)=a(Z,\1Z)\rho(Z,\1Z)$ in a
neighborhood of $p_0$ in $\C^N$, where
$a(Z,\1Z)$ is a $d\times d$ invertible matrix,
whose entries are real valued, real-analytic
functions in a neighborhood of $p_0$.

\subsection{CR points of $M$}\label{CR}
For $p\in M$ near $p_0$, the  {\em excess
codimension}
$r_1(p)$ of $M$ at $p$  is
defined by
\begin {equation}\Label{r1}
r_1(p):=d-\dim\ \span_\C\left\{\rho^j_Z(p,\1p):
1\le j\le d \right\}.
\end{equation}
Here $\rho_Z^j=(\d\rho^j/\d
Z_1,\ldots, \d\rho^j/\d
Z_N)\in\C^N$ denotes the complex gradient of
$\rho^j$ with respect to
$Z=(Z_1,\ldots,Z_N)$. It is easy to see that
$r_1(p)$ is independent of the choice of the
defining function $\rho$ and of the holomorphic
coordinates
$Z$.  A point $p_0\in M$ is called a CR {\em point}
(or $M$ is called CR at $p_0$) if the mapping
$p \mapsto r_1(p)$ is constant for $p$ in a
neighborhood of
$p_0$ in $M$. The submanifold $M$ is called CR if it
is CR at all its points and hence, by connectedness,
$r_1:=r_1(p)$ is constant on $M$. If in addition
$r_1=0$, then
$M$ is said to be {\em generic} in $\C^N$. We set
\begin{equation}\Label{V1}
V_1:=\{p\in M: M \text{ is not CR at } p \}.
\end{equation}
It is easy to see that the function $r_1(p)$
is upper-semicontinuous on $M$ and,
since $M$ is connected, the complement $M\setminus V_1$ agrees with the set of
all points in $M$, where $r_1(p)$ reaches
its minimum. The following lemma is a consequence
of the fact that $r_1(p)$ is upper-semicontinuous
for the Zariski topology on $M$ and its
proof is left to the reader.

\begin{Lem}\Label{}
The subset $V_1\subset M$ defined by {\rm (\ref{V1})}
is proper and real-analytic.
\end{Lem}

\subsection{The $(0,1)$ vector fields on $M$}
In order to define the functions $r_2(p)$ and $r_3(p)$,
we shall need the notion of $(0,1)$ vector fields on a real
submanifold $M\subset\C^N$.
For $M$ not necessarily CR and $U\subset M$ an open subset, we call a real-analytic
 vector field of the form $L=\sum_{j=1}^N
a_j(Z,\1Z)\frac{\d}{\d\1Z_j}$, with $a_j(Z,\1Z)$ 
real-analytic functions on $U$,  a {\em
$(0,1)$ vector field} on $U$ if 
\begin{equation}\Label{vfield}
(L\rho)(Z,\1Z)\equiv 0, 
\end{equation}
for any  local defining function $\rho(Z,\1Z)$ of $M$.
For $p\in M$, we denote by $\3T^{0,1}_{M,p}$  the
vector space of all germs at $p$ of  $(0,1)$ vector
fields on $M$ and by $\3T^{0,1}_{M}$ the
corresponding sheaf on
$M$ whose stalk at any $p$ is $\3T^{0,1}_{M,p}$.
It is easy to see that $\3T^{0,1}_{M}$ is independent of the choice
of $\rho(Z,\1Z)$.
Observe that $\3T^{0,1}_{M,p}$ is closed under
commutation and hence is a Lie algebra. If
$L$ is a $(0,1)$ vector field on an open set $U$
of $M$, i.e.\
$L\in \3T^{0,1}_{M}(U)$, denote by $L_p\in
\3T^{0,1}_{M,p}$  the germ of $L$ at
$p$ for $p\in U$.
If $M$ is a CR submanifold in $\C^N$,
the above definition of $(0,1)$ vector fields on $M$
coincides with the standard one
and in this case the sheaf $\3T^{0,1}_{M}$
is the sheaf of sections of a complex vector bundle
on $M$, called the CR {\em bundle} of $M$.

The following consequence of the coherence theorem of
     Oka-Cartan (see \cite{Oka} and \cite{Cartan}, Proposition~4)
will be essential for the proof that the subvariety
$V\subset M$ is real-analytic.

\begin{Lem}\Label{coherence}
Given $p_0\in M$, there exists a neighborhood
$U\subset M$ of $p_0$, an integer $m>0$ and $(0,1)$ vector fields
$L_1,\ldots,L_m\in\3T^{0,1}_{M}(U)$ such that
for any $p\in U$, any germ $\3L\in \3T^{0,1}_{M,p}$
can be written in the form
$$\3L=g_1 L_{1,p}+\cdots+g_m L_{m,p}$$
with $g_1,\ldots,g_m$ germs at $p$ of
real-analytic functions on $M$.
\end{Lem}

\begin{proof} For $p\in M$ denote by $\3A_{M,p}$
       the ring of germs at $p$ of real-analytic
functions on $M$.
For $p$ near $p_0$, we can think of
an element $L=\sum_{j=1}^N
a_j\frac{\d}{\d\1Z_j}$
 in $\3T^{0,1}_{M,p}$ as an $N$-tuple
$(a_1,\ldots,a_N)\in\3A_{M,p}^N$ satisfying
the condition in (\ref{vfield}).
Hence the subsheaf
$\3T^{0,1}_{M}\subset\3A_M^N$ coincides with the
sheaf of relations
$$\sum_{j=1}^N
a_j\Big(\frac{\d\rho^r}{\d \1Z_j}\Big)_p=0,\quad r=1,
\ldots, d.$$ Since the sheaf $\3A_M$ is coherent by
the theorem of Oka-Cartan,
it follows that $\3T^{0,1}_{M}$
is locally finitely generated
over $\3A_M$ which proves the lemma.
\end{proof}

\subsection{Degeneracy and orbit codimension}
\Label{degen}
As above let $p_0\in M$ be fixed and $\rho(Z,\1Z)$ be a local
defining function of $M$ near $p_0$. For $p\in M$
near $p_0$, we define a vector subspace $E(p)\subset
\C^N$ by
\begin{equation}\Label{rp}
E(p):=\span_\C \left\{(\3L_1\ldots
\3L_s\rho^j_Z)(p,\1p) :  1\le j\le d;
0\le s <\infty;
\3L_1,\ldots,\3L_s\in
\3T^{0,1}_{M,p} \right\}.
\end{equation}
As before $\rho_Z^j(Z,\1Z)\in\C^N$ denotes the
complex gradient of $\rho^j$ with respect to $Z$.
We leave it to the reader to check that $E(p)$ is
independent of the choice of the defining function
$\rho$ and its dimension is independent of the choice
of holomorphic coordinates $Z$ near $p$.
We call the number
\begin{equation}\Label{r2}
r_2(p):=N-\dim_\bC E(p),
\end{equation}
the {\em degeneracy} of $M$ at $p$.
We say that $M$ is of {\em
minimum degeneracy} at $p_0$
if $p_0$ is a local minimum of the function
$p\mapsto r_2(p)$.
If $r_2(p_0)=0$, we say that $M$ is
{\em finitely nondegenerate} at $p_0$.
       We
say that
$M$ is
$l$-{\em nondegenerate} at $p_0$ if $M$ is
finitely nondegenerate at $p_0$ and $l$ is the
smallest integer such that the vectors
$(\3L_1\ldots \3L_s\rho^j_Z)(p_0,\1p_0)$
span $\C^N$ for
$0\le s\le l$ and $1\le j\le d$. When
$M$ is generic, the latter definition coincides
with the one given in \cite{BER}.
(See also \cite{BHR}, where this notion appeared
for the first time for the case of a hypersurface.)

We denote by $\C T_pM:=\C\otimes_\bR T_pM$ the
complexified tangent space of $M$ at $p$ and by
$\1{\3T^{0,1}_{M,p}}$ the complex conjugates of
elements in $\3T^{0,1}_{M,p}$.
Let $\frak g_M(p)$ be the complex vector
subspace of $\C T_pM$ generated by the values at $p$
of the germs of vector fields in $\3T^{0,1}_{M,p}$,
$\1{\3T^{0,1}_{M,p}}$  and
all their commutators.
We call
\begin{equation}\Label{r3}
r_3(p):=\dim_\bR M-\dim_\C\frak g_M(p)
\end{equation}
the {\em orbit codimension} of $M$ at $p$
and say that $M$ is of {\em minimum orbit
codimension} at $p_0$ if $p_0$ is a local minimum of
the function $p\mapsto r_3(p)$.
The use of this terminology will be justified in
\S\ref{orbits}. We say that $M$ is of {\em finite
type} at $p_0$, if $r_3(p_0)=0$. When $M$ is
generic, this definition coincides with the finite
type condition of Kohn~\cite{Kohn} and
Bloom-Graham~\cite{BG}.

The following result can be obtained by applying
Lemma~\ref{coherence}, using the fact that
$\3T^{0,1}_{M,p}$ is a Lie algebra and by induction
on $s\ge 0$ in (\ref{rp}). We leave the details to the
reader.

\begin{Lem}\Label{span}
For $M\subset\C^N$, $p_0\in M$ and $\rho(Z,\1Z)$ as above, there
exist an open neighborhood
$U$ of $p_0$ in
$M$, an integer $m>0$, and
$L_1,\ldots,L_m\in
\3T^{0,1}_{M}(U)$  such that for every $p\in U$,
one has
\begin{equation}\Label{ep}
E(p)=\span_\C\left\{(L^\a\rho^j_Z)(p,\1p):\a\in\bZ_+^m;
1\le j\le d\right\},
\end{equation}
where $L^\a:=L_1^{\a_1}\ldots L_m^{\a_m}$,
$\a=(\a_1,\ldots,\a_m)$, and
\begin{multline}\Label{gp}
\frak g_M(p)=\\\span_\C\left\{
[X_{i_1},\ldots,[X_{i_{r-1}},X_{i_r}]\ldots](p):
r\ge 1;
X_{i_j}\in\{L_1,\ldots,L_m,\1L_1,\ldots,\1L_m\}
\right\}.
\end{multline}
\end{Lem}

\begin{Pro}\Label{v23}
Let $M\subset\C^N$ be a connected real-analytic
submanifold. Then the subsets $V_2, V_3\subset M$
given by
\begin{equation}\Label{V2}
V_2:=\{p\in M: M \text{ is not of minimum
degeneracy at } p \}\end{equation}
and
\begin{equation}\Label{V3}
V_3:=\{p\in M: M \text{ is not of minimum orbit
codimension at } p
\}\end{equation}
        are proper and real-analytic.
\end{Pro}

\begin{proof}
Define $r_i:=\min_{p\in M} r_i(p)$, $i=2,3$,
where $r_2(p)$ and $r_3(p)$ are the integer valued
functions defined by (\ref{r2}) and (\ref{r3}) respectively.
Given
$p_0\in M$, choose
$U$ and
$L_1,\ldots,L_m\in\3T^{0,1}_{M}(U)$
as in Lemma~\ref{span}.
Now consider the set of vector valued real-analytic
functions $L^\a\rho^j_Z$ (as in (\ref{ep})) defined
in $U$. For each subset of $N-r_2$ functions in this
set, we take all possible $(N-r_2)\times (N-r_2)$
minors extracted from their components. Then by
Lemma~\ref{span}, the set
$V_2\cap U$ is given by the vanishing of all such
minors. Since $p_0\in M$ is arbitrary,
$V_2\subset M$ is a real-analytic subvariety.
To show that $V_3\subset M$ is also real-analytic,
we repeat the above argument for the set of vector
valued real-analytic functions
$p\mapsto
[X_{i_1},\ldots,[X_{i_{r-1}},X_{i_r}]\ldots](p)$ (as
in (\ref{gp})).
Both subsets $V_2,V_3\subset M$ are proper by
the choices of $r_2$ and $r_3$.
\end{proof}

\begin{Rem}{\rm For $M\subset\C^N$ a connected
real-analytic submanifold, it follows from the
definition of
$r_1(p)$ and from the proof  of Proposition~\ref{v23}
that the sets
$\{p\in M : r_i(p)\le s\}$, $i=1,2,3$, are also
real-analytic subvarieties of $M$ for any integer $s\ge
0$. In particular, for each $i=1,2,3$, the function
$r_i(p)$ is constant in $M\setminus V_i$.}
\end{Rem}

\subsection{CR orbits in real-analytic submanifolds}
\Label{orbits}
Let $M\subset\C^N$ be a real-analytic
submanifold (not necessarily CR) and $p_0\in M$. By
a  CR {\em orbit of
$p_0$ in
$M$}  we mean a germ  at $p_0$ of a real-analytic
submanifold  $\Sigma\subset M$ through $p_0$ such
that  $\C T_p\Sigma=\frak g_M(p)$ for all $p\in
\Sigma$.  The existence (and uniqueness) of
the CR orbit of any point in $M$ follows by applying a theorem of
Nagano~(\cite{Nag}, see also \cite{BER}, \S3.1)
to  the Lie
algebra spanned by the real and imaginary parts
of the vector fields $L_1,\ldots,L_m$ given by
Lemma~\ref{span}.
The terminology introduced above for the orbit
codimension is justified by the fact that
the (real) codimension of the CR orbit of $p$ in $M$
coincides with the (complex) codimension of $\frak
g_M(p)$ in $\C T_pM$.

\section{Local structure of $M$ at a  point
of $M\setminus V$} We keep the notation
introduced in \S\ref{local}. As before we let
$r_i= \min_{p\in M} r_i(p),\  i = 1, 2, 3. $  The
following proposition gives the local structure of a
manifold near a point which is CR and also of minimum
degeneracy.

\begin{Pro}\Label{dec} Let $M\subset\bC^N$ be a
connected real-analytic submanifold and $p_0\in
M\setminus (V_1\cup V_2)$.
Then there exist local holomorphic coordinates
$Z=(Z^1,Z^2,Z^3)\in
\bC^{N_1}\times\bC^{N_2}\times\bC^{N_3}$
vanishing at $p_0$ with
$$N_1:= N -r_1-r_2,
\quad N_2:=r_2,
\quad N_3:=r_1,$$
a generic real-analytic submanifold
$M_1\subset\bC^{N_1}$ through $0$,
finitely nondegenerate at $0$,
and an open neighborhood $\3O\subset\C^N$ of $p_0$
such that
$$M\cap \3O = \big\{(Z^1,Z^2,Z^3)\in \3O:
Z^1\in M_1, Z^3=0\big\}.$$
Equivalently, in the coordinates $Z$, the germ
of
$M$ at $0$  and that of
$M_1\times\C^{N_2}\times\{0\}$ coincide.
\end{Pro}

\begin{Rem}{\rm Suppose that $M\subset\bC^N$ is
a connected real-analytic submanifold and $p_0\in
M\setminus V_2$, i.e.\ $M$ is of minimum degeneracy
at $p_0$ but not necessarily CR.  One can still ask
whether there exist local
holomorphic coordinates $Z=(Z_1,Z_2) \in \bC^{N_1}
\times \bC^{N_2}$, vanishing at $p_0$, with $N_2=r_2$
and a submanfold
$M_1 \subset \bC^{N_1}$ through $0$ which is
finitely nondegenerate at $0$ such that, in the
coordinates $Z$, the germ of
$M$ at $0$  and that of
$M_1\times\C^{N_2}$ coincide.
Observe that Proposition \ref{dec} implies that this
is the case if, in addition,
$M$ is CR at $p_0$. The following example shows that
it is not the case in general. Let
$M \subset \bC^3$ be given by $M:=\{(z_1,z_2,w)\in
\bC^3: w= z_1
\bar z_2\}$. $M$ is CR precisely at those points
where $z_1\not=0$. The $(0,1)$ vector fields on
$M$ are multiples of
$L=\partial_{\bar z_1} + z_2\partial_{\bar w}$. The
degeneracy is everywhere 1 and  the orbit dimension
is everywhere 2.  However, as the reader can easily
check,  the answer to the question above is
negative in this example with $p_0=0$.}
\end{Rem}

\begin{proof} [Proof of Proposition {\rm\ref{dec}}]
We may assume $p_0=0$. Since $M$ is CR at $0$,
there is a neighborhood of $0$ in $\bC^N$ such
that the piece of
$M$ in that neighborhood is contained as a generic
submanifold in a complex submanifold of
$\C^N$ (called the {\em intrinsic
complexification} of $M$) of
complex dimension
$N - r_1$ (see e.g. \cite{BER}).
By a suitable choice of holomorphic coordinates
$(Z^1,Z^2,Z^3)\in
\bC^{N_1}\times\bC^{N_2}\times\bC^{N_3}$
with $N_1,N_2,N_3$ as in the proposition,
we may assume that the intrinsic complexification is
given by $Z^3=0$ near $0$. Then $M$ is a generic
submanifold of the subspace $\{Z^3=0\}$.
Hence in the rest of the proof it suffices to
assume that
$M\subset\C^N$ is generic and $0\in M$.
We may therefore find holomorphic coordinates
$Z=(z,w)\in\C^n\times\C^d$, with $d$ being the
codimension of $M$ in $\C^N$ and $n:=N-d$, such
that if
$\rho=(\rho_1,\ldots,\rho_d)$ is a local defining
function of $M$ near $0$, then $\rho_w(0)$ is an
invertible $d\times d$  matrix. By the implicit
function theorem, we can write $M$ near $0$ in the
form
$$M=\{(z,w) : w-Q(z,\1z,\1w)=0\}=
\{(z,w) : \1w-\1Q(\1z,z,w)=0\},$$
where $Q$ is a $\C^d$-valued holomorphic function
defined in a neighborhood of $0$ in $\C^{2n+d}$ and
vanishing at $0$.
We now apply the definition of minimum degeneracy
given in \S\ref{degen} to the (complex valued)
defining function of $M$ given by
\begin{equation}\Label{Th}
\Theta(z,w,\1z,\1w):=\1w-\1Q(\1z,z,w).
\end{equation}
It can be easily checked that the identity
(\ref{ep}) holds with
$\rho(z,w,\1z,\1w)$ replaced by
$\Theta(z,w,\1z,\1w)$  (even though here
$\Theta(z,w,\1z,\1w)$ is complex valued). Consider
the basis of
$(0,1)$ vector fields on $M$ given by
$$L_j:=\frac{\d}{\d\1z_j} + \sum_{i=1}^d
\1Q_{\1z_j}^i(\1z,Z)\frac{\d}{\d\1w_i},
\quad 1\le j\le n,$$
where, as above, $Z=(z,w)$.
Observe that since $\1Q$ is independent
of $\1w$, for
$\a\in\bZ_+^n$ and
$1\le j\le d$,
$$L^\a\Theta^j_Z(Z,\1Z)=-\1Q^j_{Z,\1z^\a}(\1z,Z).$$
Since $M$ is
of minimum degeneracy at $0$, it follows  that
for $Z$ in a neighborhood of $0$ in $M$,
$$\dim\, \span_\C
\{\1Q^j_{Z,\1z^\a}(\1z,Z) : \a\in\bZ_+^n; 1\le j\le
d \} = N-r_2.$$
By a standard complexification argument (see e.g.
Lemma 11.5.8 in \cite{BER}), we conclude that for
$\z\in\C^n$ and
$Z\in\C^N$ near the origin, we also have
$$\dim\, \span_\C
\{\1Q^j_{Z,\z^\a}(\z,Z) : \a\in\bZ_+^n; 1\le j\le
d \} = N-r_2.$$
Hence there exists an integer $l\ge 0$ such that
for $Z\in\C^N$ in a neighborhood of $0$,
$$\dim\, \span_\C
\{\1Q^j_{Z,\z^\a}(0,Z) : 0\le|\a|\le l; 1\le
j\le d \} = N-r_2.$$
In particular, if $K$ is $d$ times the number of
multi-indices $\alpha \in \bZ_+^n$ with $0\le |\a|
\le l$, the map
$\psi$ given by
$$Z\mapsto \psi(Z)=
(\1Q^j_{\z^\a}(0,Z))_{0\le|\a|\le l, 1\le j\le
d }\in \C^K$$ is of constant rank equal to
$N-r_2$ for $Z$ in a neighborhood of $0$ in
$\C^N$. By the implicit function theorem,
there exists a holomorphic change of
coordinates
$Z=\Phi(\2Z^1,\2Z^2)$ with
$(\2Z^1,\2Z^2)\in\C^{N-r_2}\times\C^{r_2}$
such that $\psi(\Phi(\2Z^1,\2Z^2))\equiv
\psi(\Phi(\2Z^1,0))$.
It follows that $\1Q_{\z^\a}(0,\Phi(\2Z^1,\2Z^2))$
is independent of $\2Z^2$ for all $\a$,
$0\le|\a|\le l$, and hence, by the choice of
$l$, for all $\a\in\bZ_+^n$. Therefore, if we
write the complexification of (\ref{Th}) in
the form
$$\Theta(Z, \Z)= \tau - \1Q(\chi,z,w),\ \
Z=(z,w)\in
\bC^n\times \bC^d,\
\Z = (\chi, \tau) \in \bC^n\times \bC^d,$$ we
conclude  that
$\Theta(\Phi(\2Z^1,\2Z^2),\Z)$ is independent of
$\2Z^2$. Hence the (complex valued) function given
by
$$\2\Theta(\2Z^1,\2Z^2,\2\Z^1,\2\Z^2):=
\Theta(\Phi(\2Z^1,\2Z^2),\1{\Phi}(\2\Z^1,\2\Z^2))$$
is independent of $\2Z^2$, and $M$ is given by
$\2\Theta\big(\2Z^1,\2Z^2,\1{\2Z^1},\1{\2Z^2}
\big)=0$. Thus all vector fields $\d/\d\2Z_j^2$,
$1\le j\le r_2$, are tangent to $M$ and
hence so are the vector fields $\d/\d\1{\2Z^2_j}$.
After a linear change of the coordinates
$\2Z^1=(\2Z^{11},\2Z^{12})\in\C^{n-r_2}\times\C^{d}$
we can write
$M$ near $0$ in the form
$$M=\{ (\2Z^{11},\2Z^{12},\2Z^2) :
\Im\ \2Z^{12}=\phi(\2Z^{11},\1{\2Z^{11}},\Re\
\2Z^{12})
\},
$$
where $\phi$ is a real-analytic, real
vector valued function. Hence the submanifold
$M_1\subset\C^{N_1}$ given by
$M_1:=M\cap\{\2Z^2=0\}$ satisfies the required
assumptions.
\end{proof}

The following proposition gives the structure of a
generic submanifold at a point of minimum orbit
codimension.  Recall that we have used the notation
$r_3 =
\min_{p\in M} r_3(p)$, where $r_3(p)$ is the orbit
codimension of
$p$.

\begin{Pro}\Label{ff}
Let $M\subset \bC^N$ be a be a
connected real-analytic generic submanifold and
$p_0\in M$. The following are equivalent:
\begin{enumerate}
\item [(i)] $p_0 \in M\setminus V_3$
\item [(ii)] There is an open neighborhood $U$ of
$p_0$ in $M$ and a real-analytic mapping
$$h: U\to \bR^{r_3}, \ \ h(p_0) = 0,
$$
of rank $r_3$, which extends holomorphically to
an open neighborhood of $U$ in $\bC^N$, such
that
$h^{-1}(0)$ is a CR manifold of finite type.
\item [(iii)] In addition to the assumptions
of the condition \hbox {\rm (ii)},  for all
$u$ in a neighborhood of
$0$ in $\bR^3$,
$h^{-1}(u)$ is a CR manifold of finite type.
\end{enumerate}
\end{Pro}
\begin{proof}
Since $M$ is generic, and hence CR, we
can choose a frame $(L_1,\ldots,L_n)$
       of
real-analytic
$(0,1)$ vector fields on $M$ near $p_0$, spanning
the space of all $(0,1)$ tangent vectors to $M$ at
every point near $p_0$. (Here $n=N-d$, where
$d$ is the codimension of $M$.) We write
$L_j=X_j+\sqrt{-1}X_{j+n}$, where $X_j$,
$1\le j\le 2n$, are real valued vector fields. We
prove first that (i) implies (iii).   By the
condition that
$M$ is of minimum orbit codimension $r_3$ at
$p_0$, it follows that the collection of the
vector fields
$X_j$,
$1\le j\le 2n$, generates a Lie algebra, whose
dimension at every point  near
$p_0$ is
$2n+d -r_3$. Therefore, by the (real) Frobenius
theorem, we conclude that there exist $r_3$
real-analytic real valued functions
$h_1,\ldots,h_{r_3}$ with independent
differentials, defined in a neighborhood of
$p_0$, vanishing at $p_0$ and such that
$L_jh_m\equiv 0$ (i.e. $h_m$ is a CR function) for all $1\le j\le n$ and
$1\le m\le r_3$. Moreover, the local orbits of
the $X_j$, $1\le j\le 2n$, are all of the form
$M_u=\{p\in M: h(p)=u\}$ with
$h=(h_1,\ldots,h_{r_3})$ and $u\in\bR^{r_3}$
sufficiently small. By a theorem of Tomassini 
(\cite{Tomassini}, see also \cite{BER}, Corollary~1.7.13), the
functions $h_1,\ldots,h_{r_3}$ extend
holomorphically to a full neighborhood of
$p_0$ in $\C^N$. This proves that (i) implies
(iii). For the proof that (ii) implies (i),
we observe that since  $h$ extends
holomorphically, we have $L_j h_m \equiv 0$
for all $1\le j\le n$ and
$1\le m\le r_3$.  By the reality of $h_m$ it
follows that
$X_j h_m \equiv 0$ for all $1\le j\le 2n$ and
$1\le m\le r_3$. Hence the  set $M_0:= h^{-1}
(0)$ is the CR orbit of $M$ at $p_0$ and is
of dimension
$r_3$, which proves (i). Since the
implication (iii) $\implies$ (ii) is trivial,
the proof of the proposition is complete.
\end{proof}

The following proposition gives useful local
holomorphic coordinates for a generic
submanifold around a point of minimum orbit
codimension.

\begin{Pro}\Label{dec-orbit} Let
$M$ be a connected generic real-analytic submanifold
of $\C^N$ of codimension $d$ and $p_0\in
M\setminus V_3$. Set $n:=N-d$, $d_2:=r_3$ and
$d_1:=d-r_3$. Then there exist
holomorphic local coordinates
$Z=(z,w,u)\in \C^N=\bC^{n}\times\bC^{d_1}\times\bC^{d_2}$
vanishing at $p_0$, an open neighborhood
$\3O=\3O_1\times\3O_2\subset\bC^{n+d_1}\times\bC^{d_2}$
of $p_0$, and a holomorphic map $Q$ from a
neighborhood of $0$ in
$\C^n\times\C^n\times\C^{d_1}\times\C^{d_2}$ to
$\bC^{d_1}$ satisfying
\begin{equation}\Label{Qcond} Q(z,0,\tau,u)\equiv
Q(0,\chi,\tau,u)\equiv\tau
\end{equation}
such that
$$M\cap\3O=\big\{(z,w,u)\in\3O: u\in\bR^{d_2},
\, w=Q(z,\1z,\1w,u)\big\},$$
and for every
$u\in\bR^{d_2}$ close to $0$ the submanifold
\begin{equation}\Label{mu}
M_u:=\big\{(z,w)\in \3O_1:w=Q(z,\bar z,\bar
w,u)\big\}\subset\C^{n+d_1}
\end{equation}
is generic and of finite type.
\end{Pro}

\begin{proof}
We take normal coordinates $Z'=(z',w')\in
\C^n\times\C^d$ vanishing at $p_0$ (see e.g. \cite{BER}, \S4.2), i.e. we assume
that $M$ is given by $w'=Q'(z',\1z',\1w')$
near
$0$, where
$Q'$ is a germ at $0$ in
$\bC^{2n+d}$ of a holomorphic
$\bC^{d}$-valued function satisfying
\begin{equation}\Label{Qcond1}
Q'(z',0,\tau')\equiv
Q'(0,\chi',\tau')\equiv\tau'.
\end{equation}
We may choose a frame $(L_1,\ldots, L_n)$
spanning the space of all $(0,1)$ vector
fields on $M$ of the form
$L_j=\frac{\d}{\d\1z'_j}+\sum_{i=1}^d
{Q'}^i_{\1z'_j}(z',\1z',\1w')\frac{\d}{\d\1w'}$
for
$1\le j\le n$. In particular, $L_j(0)=
\frac{\d}{\d\1z'_j}$. 
Let $h= (h_1,\ldots, h_{d_2})$ be the
functions given by (iii) in
Proposition \ref{ff}.
Since, for
$1\le m\le d_2$, the  functions
$h_m$ are real and extend holomorphically, we
conclude that
$L_jh_m\equiv
\1L_jh_m\equiv 0$.  We denote again the by
$h_1,\ldots, h_{d_2}$ the extended functions. By
the choice of the coordinates, $\d h_m/\d
z'_j(0)=0$, $1\le m\le d_2$, $1\le j\le n$. By
using the independence of the differentials of
$h_1,\ldots,h_{d_2}$ and reordering
the components $w'_1,\ldots,w'_{d}$ if
necessary, we may assume that
$$\det \Big(\frac{\d h_m}{\d w'_j}(0)\Big)_{1\le
m\le d_2,\, d_1+1\le j \le d}\ne 0.$$ We make the
following change of holomorphic coordinates in
$\C^N$ near $0$:
\begin{equation}\Label{}
z''=z', \quad w''_j=w'_j \text{ for } 1\le j\le
d_1,
\quad  w''_j=h_{j-d_1}(z',w') \text{ for }
d_1+1\le j\le d.
\end{equation}
Note that on $M$, we have $w_j''=\1w_j''$ for
$d_1+1\le j\le d$. The reader can check that the
new coordinates
$(z'',w'')\in\C^n\times\C^d$ are again normal for
$M$. Indeed,
$M$ is given by $w''=Q''(z'',\1z'',\1w'')$ where
$Q''$ satisfies  the analog of (\ref{Qcond1}),
with
${Q''}^j(z'',\1z'',\1w'')\equiv\1w''_j$ for
$d_1+1\le j\le d$. The desired coordinates are
obtained by taking
$(z,w,u):=(z'',w'')$ i.e.\ $z=z''$ and $(w,u) =
w''$ with $w\in \bC^{d_1}$,
$u\in\C^{d_2}$. We  take $Q^j: = {Q''}^j$ for $1
\le j \le d_1$. By the properties of the
functions
$h_1,\ldots,h_{d_2}$, the submanifold
$M_u$ given by (\ref{mu}), with $u\in\bR^{d_2}$ close to $0$,
is of finite type if $\3O_1$ is a sufficiently
small neighborhood of $0$ in $\bC^{n+d_1}$. This
completes the proof of Proposition
\ref{dec-orbit}.
\end{proof}


\section{Properties of $k$-equivalences between
germs of real submanifolds}\Label{kequiv}

We first observe that if $(M,p)$ and $(M',p')$
are two germs in $\bC^N$ of real-analytic
submanifolds at
$p$ and $p'$ respectively, then for any
formal
$k$-equivalence
$H$ between
$(M,p)$ and $(M',p')$, the $k$th Taylor
polynomial of $H$ is a convergent
$k$-equivalence. Therefore we may, and shall,
assume all
$k$-equivalences in the rest of this paper to be convergent.
By a {\em local parametrization} $Z(x)$ of $M$
at $p$ we shall mean a real-analytic
diffeomorphism $x\mapsto Z(x)$ between open
neighborhoods of $0$ in $\bR^{\dim M}$  and of
$p$ in $M$ satisfying $Z(0)=p$.  We say that a
function $f(x)$ in a neighborhood of $0$ in
$\bR^m$
       {\em vanishes of order} $k$ at $0$, if $f(x)=O(|x|^k)$.

One of the main results of this section is to
show that $k$-equivalences, for sufficiently
large $k$ preserve the integers $r_j(p)$, $j=1,
2, 3$ introduced in \S2, and their minimality.
For simplicity of notation we state the result
for $p=p'=0$.

\begin{Pro}\Label{const}
Let $(M,0)$
and $(M',0)$ be two germs  at $0$ in
$\C^N$ of
real-analytic submanifolds which are
$k$-equivalent for every $k$. Denote by
$r_j(0)$ and
$r'_j(0)$,
$j=1,2,3$, the integers given by (\ref{r1}),
(\ref{r2}), and
(\ref{r3}) for $M$ and $M'$ respectively.
Then the following hold.
\begin {enumerate}
\item  [(i)] $r_1(0)=
r'_1(0)$. Also $M$ is CR at $0$
if and only if $M'$ is CR at $0$.
\item [(ii)]If $M$ is CR at $0$ then $r_2(0)=
r'_2(0)$, and
$M$ is of minimum degeneracy at $0$ if and
only  if $M'$ is of  minimum degeneracy at $0$.
\item [(iii)]If $M$ is CR at $0$ then $r_3(0)=
r'_3(0)$, and
$M$ is of minimum orbit codimension at $0$ if and
only  if $M'$ is of  minimum orbit codimension
at $0$.
\end {enumerate}
\end{Pro}

Before proving Proposition \ref{const}, we shall
need some preliminary results. The following
useful but elementary lemma gives alternative
definitions of
$k$-equivalences.

\begin{Lem}\Label{equivss}
Let $H\colon(\C^N,0)\to(\C^N,0)$ be an invertible
germ of a holomorphic map and $(M,0)$ and $(M',0)$
be
two germs at $0$ of real-analytic submanifolds of
$\C^N$ of the same dimension. Then for any
integer
$k>1$, the following are equivalent:
\begin{enumerate}
\item[(i)]
$H$ is a $k$-equivalence between $(M,0)$ and
$(M',0)$.
\item[(ii)]
There exist local
parametrizations $Z(x)$ and $Z'(x)$ at $0$
of $M$ and $M'$ respectively such that
$Z'(x)=H(Z(x))+O(|x|^k)$.
\item[(iii)]
For every local
parametrization $Z(x)$ of $M$ at $0$,
there exists a local parametrization $Z'(x)$
of $M'$ at $0$ such that $Z'(x)=H(Z(x))+O(|x|^k)$.
\item[(iv)]
There exist local defining functions $\rho(Z,\1Z)$ and
$\rho'(Z',\1{Z'})$ of $M$ and $M'$ respectively near $0$ such that
$\rho'(H(Z),\1H(\Z))=\rho(Z,\Z)+ O(|(Z,\Z)|^k)$.
\item[(v)]
For every local defining function $\rho(Z,\1Z)$ of $M$ near $0$,
there exists a local defining function $\rho'(Z',\1{Z'})$ of
$M'$ near $0$ such that
$\rho'(H(Z),\1H(\Z))=\rho(Z,\Z)+ O(|(Z,\Z)|^k)$.
\item[(vi)]
For any local defining functions $\rho(Z,\1Z)$
and
$\rho'(Z',\1{Z'})$ of $M$ and $M'$ respectively near $0$, there exists a
holomorphic function $a(Z,\Z)$  defined in a
neighborhood of $0$ in $\C^{2N}$ with values
in the space of $d\times d$ invertible matrices
(where $d$ is the codimension of $M$) such that
$\rho'(H(Z),\1H(\Z))=a(Z,\Z)\rho(Z,\Z)+
O(|(Z,\Z)|^k)$.
\end{enumerate}
In particular,
inverses and compositions of
$k$-equivalences are also $k$-equivalences.
\end{Lem}

Since the proof of Lemma \ref{equivss} is
elementary, it is left to the reader. We shall
also need the following two lemmas for the proof
of Proposition
\ref{const}.

\begin{Lem}\Label{va}
Let $(v_\a(x))_{\a\in A}$  be a
collection of real-analytic
$\bC^K$-valued functions in a neighborhood of $0$ in
$\bR^m$. If the dimension of the span in
$\bC^K$ of the $v_\a(x)$, $\a\in A$, is not constant for $x$ in
any neighborhood of $0$, then there exists an integer
$\k>1$ such that, for any other collection
of real-analytic $\bC^K$-valued functions
$(v'_\a(x))_{\a\in A}$ in some neighborhood of $0$
with
$v'_\a(x)=v_\a(x)+O(|x|^\k)$, the dimension of the
span in
$\bC^K$ of the
$v'_\a(x)$ is also nonconstant in any neighborhood
of $0$.
\end{Lem}

\begin{proof}
Denote by $r$ the dimension of the span in
$\bC^K$ of the
$v_\a(0)$, $\a\in A$. By the assumption, there
exists an
$(r+1)\times (r+1)$ minor $\Delta(x)$, extracted
from the components of the $v_\a(x)$, which does not
vanish identically. Note that $\Delta(0)=0$. Let
$\gamma\in\bZ_+^m$, $|\gamma|\ge 1$, be  such
that
$\d^\gamma\Delta(0)\ne 0$. Then for $\k:=|\gamma|+1$
and $v'_\a(x)$ as in the lemma, it follows that
$\d^\gamma\Delta'(0)\ne 0$, where $\Delta'(x)$ is
the corresponding minor with $v_\a(x)$ replaced by
$v'_\a(x)$. On the other hand, the dimension of the
span of the $v'_\a(0)$ is also $r$. Since
$\Delta'(x)$ is an $(r+1)\times (r+1)$ minor
that does not vanish identically, the proof of the
lemma is complete.
\end{proof}

\begin{Lem}\Label{add0}
Let $M_1,M'_1\subset\C^{N_1}$ be two generic
real-analytic submanifolds
through
$0$, of the same dimension. Let
$M:=M_1\times\{0\}$ and
$M':=M'_1\times\{0\}$, both contained
in $\C^N=\C^{N_1}\times\C^{N_2}$,
and  $H$  a
$k$-equivalence between
$(M,0)$ and $(M',0)$,
with $k>1$.
        Let $Z=(Z^1,Z^2)$ and $H=(H^1,H^2)$ be
the corresponding decompositions for the components
of $Z$ and $H$. Then $H^2(Z^1,0)=O(|(Z^1)|^k)$ and
$Z^1\mapsto H^1(Z^1,0)$ is a $k$-equivalence
between $(M_1,0)$ and $(M'_1,0)$.
\end{Lem}

\begin{proof}
We write $Z'=(Z'{}^1,Z'{}^2)\in
\C^{N_1}\times\C^{N_2}$. Let
$\rho'_1(Z'{}^1,\1{Z'{}^1})$ be a local defining
function for $M'_1\subset\C^{N_1}$. Then
\begin{equation}\Label{rho'-1}
\rho'(Z',\1{Z'}):=(\rho'_1(Z'{}^1,\1{Z'{}^1}),
\Re\, Z'{}^2, \Im\, Z'{}^2)
\end{equation}
is a local defining function for $M'$ in a
neighborhood of $0$ in $\C^N$.
By the definition of $k$-equivalence, we obtain
\begin{equation}\Label{kids-1}
\rho'_1(H^1(Z(x)),\1{H^1(Z(x))})=O(|x|^k),
\quad H^2(Z(x))=O(|x|^k),
\end{equation}
for any local parametrization $Z(x)$ of
$M$ at $0$. By the second identity in (\ref{kids-1}), the
holomorphic function $H^2(Z^1,0)$ vanishes of
order
$k$ at
$0$ on the submanifold $M_1$ which
is generic in  $\C^{N_1}$.
This implies the first statement of the lemma.
       Since $H$ is invertible and by the first
statement of the lemma we have $H^2_{Z}(0)=0$,
the map
$Z^1\mapsto H^1(Z^1,0)$ must be invertible at
$0$. Hence the first identity in (\ref{kids-1})
implies the second statement of the lemma.
\end{proof}

\begin{proof} [Proof of
Proposition~{\rm\ref{const}}]We first observe
       that every
$2$-equivalence between $(M,0)$ and $(M',0)$
induces a linear isomorphism between
$T_0M$ and $T_0M'$. Since $(M,0)$ and $(M',0)$
are $k$-equivalent for every $k$, this implies
$r_1(0)=r_1'(0)$. To complete the proof of (i),
we argue by contradiction. We assume that $M'$ is
CR at
$0$ but that $M$ is not. If $\rho(Z,\1Z)$ is a
local defining function for $M$ and $Z(x)$ is a
local parametrization of $M$ at $0$, we set
$v^j(x):=\rho^j_Z(Z(x),\1{Z(x)})$, $1\le j\le d$.
Since $M$ is assumed not to be CR at $0$, the
collection of functions
$v^j(x)$ satisfies the assumptions of
Lemma~\ref{va}. Let $\k$ be the integer given by
the lemma. We take $k\ge\k+1$ and let $H$  be a
$k$-equivalence between $(M,0)$ and $(M',0)$.
If we set $\2M:=H^{-1}(M')$, then the identity map
is a $k$-equivalence between $(M,0)$ and
$(\2M,0)$. Hence, by Lemma~\ref{equivss} (iii,v),
there exist a local parametrization $\2Z(x)$
of $\2M$ at $0$ and a local defining function
$\2\rho$ for $\2M$ near $0$ such that
$\2Z(x)=Z(x)+O(|x|^k)$ and
$\2\rho(Z,\1Z)=\rho(Z,\1Z)+O(|Z|^k)$.
We apply Lemma~\ref{va} for the collection
$v^j(x)$ defined above and
$v'{}^j(x):=\2\rho^j_Z(\2Z(x),\1{\2Z(x)})$ and
conclude that $\2M$ is not CR at $0$.
Thus we have reached a contradiction, since $\2M$
and
$M'$ are biholomorphically equivalent.
This completes the proof of (i).

To prove (ii), suppose that $M$
and $M'$ are CR at $0$. Since $M$ and $M'$ are CR and $r_1(0)= r_1'(0)$
by (i), we may assume that
$M=M_1\times\{0\}$ and $M'=M'_1\times\{0\}$,
both contained in $\C^{N_1}\times\C^{N_2}$ with
$N_1:=N -r_1(0)$, $N_2:=r_1(0)$ and $M_1$ and
$M'_1$ generic in $\bC^{N_1}$ (cf.\ beginning of
proof of Proposition~\ref{dec}). By
Lemma~\ref{add0},
$M_1$ and $M'_1$ are also $k$-equivalent for
every
$k>1$. We observe that $M$ is of minimum
degeneracy at $0$ if and only if $M_1$ is of
minimum degeneracy at $0$ and the degeneracies
of $M$ and $M_1$ at $0$ are the same. Therefore,
by replacing
$M$ by
$M_1$ and
$M'$ by
$M'_1$,  we may assume that
$M$ and $M'$ are generic (i.e. $
r_1(0)=r'_1(0)=0$) in the rest of the proof.

       We show first that
$r_2(0) = r'_2(0)$. From the definition (\ref
{r2}) of these numbers, there exists $l\ge 0$
such that
\begin{equation}\Label{el}
\dim_\bC E_l(0) = r_2(0), \ \ \
\dim_\bC E'_l(0) = r'_2(0),
\end{equation}
where, for $p \in M$,
$$E_l(p):=\span_\C\left\{(\3L_1\ldots
\3L_s\rho^j_Z)(p,\1p) :  0\le s\le l;
\3L_1,\ldots,\3L_s\in
\3T^{0,1}_{M,p}; 1\le j\le d\right\},$$
with $\rho(Z,\1Z)$ being a defining function for
$M$ near $0$, and $E'_l(p')\subset\C^N$ is the
corresponding subspace for
$M'$. We may choose holomorphic
coordinates
$Z=(z,w)\in\C^n\times\C^d$ vanishing at
$0$ such that the $d\times d$ matrix $\rho_w(0)$ is
invertible. In these coordinates we take a basis of
$(0,1)$ vector fields on $M$ in the form
\begin{equation}\Label{Lj}
        L_j=
\frac{\d}{\d\1z_j}-^\tau\!\!\rho_{\1z_j}\:
^\tau\!(\rho_w^{-1})\Big(\frac{\d}{\d\1w}\Big),
\quad 1\le j\le n,
\end{equation}
where we have used matrix notation so that
$\left(\frac{\d}{\d\1w}\right)=
\big(\frac{\d}{\d\1w_1},\ldots,\frac{\d}{\d\1w_1}
\big)$ is viewed as a $d\times 1$ matrix.
We now choose a local parametrization
$Z(x)$ of $M$ at $0$,
and put
\begin{equation}\Label{vja}
v^j_\a(x):=L^\a\rho^j_Z(Z(x),\1{Z(x)}),
\quad 1\le j\le d,\, \ \a\in\bZ_+^n,\,\ 0\le
|\a|\le l.
\end{equation}
We then choose $k>l+1$ and $H$ to be a
$k$-equivalence between $(M,0)$ and $(M',0)$ which
exists by the assumptions of the proposition.
Replacing $M'$ by $H^{-1}(M')$ we may assume without
loss of generality that $H$ is the identity map of
$\C^N$.
By Lemma~\ref{equivss} (ii,v), we can find a local
parametrization $Z'(x)$ of $M'$ at $0$ satisfying
$Z'(x)=Z(x)+O(|x|^k)$ and a local defining function
$\rho'$ of $M'$ near $0$ satisfying
$\rho'(Z,\1Z)=\rho(Z,\1Z)+O(|Z|^k)$.
Denote by $L'_j$, $1\le j\le n$, the local basis of
$(0,1)$ vector fields on $M'$ given by the analog of
(\ref{Lj}) with $\rho$ replaced by $\rho'$.
(Observe that $\rho'_w(0)$ coincides with $\rho_w(0)$
and hence is invertible). By the choice of $\rho'$,
we have
$L'_j=L_j+R_j$ in a neighborhood of $0$ in $\C^N$, where $R_j$ is a
vector field whose
coefficients vanish of order $k-1$ at $0$.
We put
\begin{equation}\Label{vja'}
v'{}^j_\a(x):=L'{}^\a\rho'{}^j_Z(Z'(x),\1{Z'(x)}),
\quad 1\le j\le d,\  \a\in\bZ_+^n, \ 0\le |\a|\le
l.
\end{equation}
Then it follows from the construction that
\begin{equation}\Label{veq}
v'{}^j_\a(x)=v^j_\a(x)+O(k-l-1)
\end{equation}
       and, in
particular, $v'{}^j_\a(0)=v^j_\a(0)$ for all
$j$ and $\alpha$ as in (\ref{vja}).  Hence, by
making use of (\ref{el}), we have
$r_2(0) = r'_2(0)$, which proves the first part
of (ii).

To prove the second part of (ii), assume that
       $M'$ is of
minimum degeneracy at $0$ and that $M$ is not.
We shall reach a contradiction by again
making use of Lemma~\ref{va}.
      From the definition of minimum degeneracy
there exists an integer
$l'\ge 0$ such that
\begin{equation}\Label{elp} \dim
E'_{l'}(p')\equiv
\dim E'_{l'}(0) , \ \ \
\dim
E_{l'}(p)\not\equiv
\dim E_{l'}(0),
\end{equation}
        for $p\in M$, and $p'\in M'$ near $0$.
Hence the collection of real-analytic
functions given by (\ref{vja}) with $l$ replaced
by $l'$ satisfies the assumption of
Lemma~\ref{va}. We let $\k>1$ be the integer
given by that lemma and choose $H$ to be a
$k$-equivalence with $k$ satisfying
$\k=k-l'-1$.  As before, we may assume that $H$
is the identity.   Using
again Lemma \ref{equivss} (ii,v), we obtain the
analogue of (\ref{veq}), with
$l$ replaced by $l'$.
        We conclude by
Lemma~\ref{va} that the dimension of the span of
the
$v'{}^j_\a(x)$ given by (\ref{vja'}), with $l$
replaced by $l'$, is not
constant in any neighborhood of $0$. This
contradicts the first part of
(\ref{elp}) and proves the second part of
(ii).

The proof of (iii) is quite similar to that of
(ii), and the details are left to the reader.
The proof of Proposition
\ref{const} is complete.
\end{proof}

\section{Reduction of Theorem~\ref{2} to the
case of generic, finitely nondegenerate
submanifolds}\Label{red}

In this section we reduce Theorem~\ref{2} to the
case where $M$ and $M'$ are generic and
$M'$ is finitely nondegenerate. For this case,
the precise statement is given in
       Theorem~\ref{tech} below. After the
reduction to this case, the rest of the paper
will be devoted to the proof of
Theorem~\ref{tech}.

\begin{Thm}[Main Technical Theorem]\Label{tech} Let
$(M,0)$ and $(M',0)$ be two germs
of generic real-analytic submanifolds of $\C^N$
of the same dimension.
Assume that $M$ is of
minimum orbit codimension  at $0$ and that
$M'$ is finitely nondegenerate at $0$. Then for
any integer
$\k>1$, there exists an integer $k>1$ such that
if $H$ is a
$k$-equivalence between $(M,0)$ and $(M',0)$, then
there exists a biholomorphic equivalence $\4H$ between $(M,0)$ and $(M',0)$
with
$\4H(Z)=H(Z) + O(|Z|^\k)$.
\end{Thm}

\begin{Rem}\Label {hyp}{\rm We should mention here
that the proof of Theorem~\ref{tech} is simpler
when
$M$ and
$M'$ are hypersurfaces of $\bC^N$.  In fact in
this case under the assumptions of the theorem,
if $H$ is a formal equivalence between $(M,0)$
and $(M',0)$, then if follows from Theorem 5 in
\cite{Asian} that $H$ is convergent.  Hence, in
particular, the proof of the equivalence of (ii)
and (iii) in Corollary \ref{main} is simpler in
the case of hypersurfaces.}
\end{Rem}

In order to show that Theorem \ref{2} is a
consequence of Theorem \ref{tech}, we shall
need the following.
\begin{Lem}\Label{add}
  Let $M_1,M'_1\subset\C^{N_1}$ be generic
real-analytic submanifolds through $0$
and assume that $M'_1$ is $l$-nondegenerate at $0$.
Let
$M:=M_1\times\C^{N_2}$ and
$M':=M'_1\times\C^{N_2}$ (both contained
in $\C^N=\C^{N_1}\times\C^{N_2}$), and let $H$ be
a
$k$-equivalence between
$(M,0)$ and $(M',0)$ with $k> l+1$.
        Let $Z=(Z^1,Z^2)$ and $H=(H^1,H^2)$ be
the corresponding decompositions for the components
of $Z$ and $H$. Then the following hold:
\begin{enumerate}
\item [(i)] $\left(\partial H^1/\partial
Z^2\right) (Z) =O(|(Z)|^{k-l-1})$;
\item [(ii)] $Z^1\mapsto
H^1(Z^1,0)$ is a
$k$-equivalence between $(M_1,0)$ and $(M'_1,0)$.
\end{enumerate}
\end{Lem}

\begin{proof} Observe that $M$ and
$M'$ are generic submanifolds of $\bC^N$.
        We write
$Z'=(Z'{}^1,Z'{}^2)\in
\C^{N_1}\times\C^{N_2}$. Let
$\rho'_1(Z'{}^1,\1{Z'{}^1})$ be a local defining
function for $M'_1\subset\C^{N_1}$. Then
$\rho'(Z',\1{Z'}):=\rho'_1(Z'{}^1,\1{Z'{}^1})$
is a local defining function for $M'$ in a
neighborhood of $0$ in $\C^N$.
By the definition of $k$-equivalence, we obtain
\begin{equation}\Label{kids}
\rho'_1(H^1(Z(x)),\1{H^1(Z(x))})=O(|x|^k).
\end{equation}
for any local parametrization $Z(x)$ of
$M$ at $0$. We choose $Z(x)$ in the form
\begin{equation}\Label{zx}
\bR^{\dim M_1}\times \bR^{
N_2}\times \bR^{N_2}\ni
x=(x^1,x^2,y^2)\mapsto Z(x)=
\big(Z^1(x^1),x^2+iy^2\big)\in M,
\end{equation}
where $x \mapsto Z^1(x^1)$ is a local
parametrization of $M_1$ at $0$.
Similarly,
we choose a local defining function
$\rho_1(Z^1,\1{Z^1})$ for $M_1$ near $0$ in
$\C^{N_1}$ and put
$\rho(Z,\1Z):=\rho_1(Z{}^1,\1{Z{}^1})$.
Since $H$ is a $k$-equivalence between $(M,0)$ and
$(M',0)$, the identity map is a $k$-equivalence
between $(M,0)$ and $(\2M,0)$ with
$\2M:=H^{-1}(M')$. We let $L_j$, $1\le j\le n$, be
the $(0,1)$ vector fields defined in a
neighborhood of $0$ in
$\C^N$ given by (\ref{Lj}) (after reordering
coordinates in the form
$Z=(z,w)\in\C^N$ with $\rho_w(0)$ invertible).
Similarly we define
$\2L_j$ by an analogue of (\ref{Lj}) with $\rho$
replaced by $\2\rho$, where $\2\rho$ is the
defining function of $\2M$ given by
Lemma~\ref{equivss} (v) for the identity map so
that
$\2\rho(Z,\1Z)=\rho(Z,\1Z)+O(|Z|^k)$. (We may take
the same decomposition $Z=(z,w)$ since
$\2\rho_w(0)=\rho_w(0)$). Hence
$\2L_j=L_j+R_j$ with
$R_j$ a
$(0,1)$ vector field in a neighborhood of $0$ in
$\C^N$
         whose coefficients vanish of order $k-1$ at $0$.
Observe that the vector fields $L'_j:=H_*
\2L_j$ are tangent to $M'$.

By Lemma~\ref{equivss} (vi), there
exists a $d\times d$ real-analytic matrix valued
function
$a(Z,\1Z)$ such that
\begin{equation}\Label{kid}
\rho'(H(Z),\1{H(Z)})=
a(Z,\1Z)\rho(Z,\1Z)+ O(|Z|^k).
\end{equation}
Differentiating (\ref{kid}) with respect to $Z$ and
applying $L^\a$ for $|\a|\le l$, we obtain
\begin{multline}\Label{kid1}
L^\a\big(\rho'_{Z'}(H(Z),\1{H(Z)}) H_Z(Z)\big)=
a(Z,\1Z) (L^\a\rho_Z)(Z,\1Z) + \\
\sum_{0\le|\b|<|\a|}
A_\b(Z,\1Z) (L^\b\rho_Z)(Z,\1Z) +
L^\a\bigg(\sum_{j=1}^d
\rho^j(Z,\1Z)B_j(Z,\1Z)\bigg) + O(|Z|^{k-l-1}),
\end{multline}
where $A_\b(Z,\1Z)$ and $B_j(Z,\1Z)$ are
real-analytic functions in a neighborhood of
$0$ in $\C^N$, valued in $d\times d$ and in $d\times
N$ matrices respectively. Using the relation $\2L_j=L_j+O(k-1)$
and the definition of $L'_j$ given above, we
conclude
\begin{multline}\Label{kid2}
(L'{}^\a\rho'_{Z'})(H(Z),\1{H(Z)}) H_Z(Z)=
a(Z,\1Z) (L^\a\rho_Z)(Z,\1Z) + \\
\sum_{0\le|\b|<|\a|}
A_\b(Z,\1Z) (L^\b\rho_Z)(Z,\1Z) +
L^\a \bigg( \sum_{j=1}^d
\rho^j(Z,\1Z)B_j(Z,\1Z) \bigg) +
O(|Z|^{k-l-1}).
\end{multline}

We now choose a local parametrization $Z'(x)$ of $M'$ at $0$ given by
Lemma~\ref{equivss} (iii), i.e.
$Z'(x)=H(Z(x))+O(|x|^k)$, with $Z(x)$ given by
(\ref{zx}). Since the $L_j$ are tangent to $M$, we
conclude from (\ref{kid2})
that
\begin{multline}\Label{kid3}
(L'{}^\a\rho'_{Z'})\big(Z'(x),\1{Z'(x)}\big) H_Z(Z(x))=
a(Z(x),\1{Z(x)}) (L^\a\rho_Z)(Z(x),\1{Z(x)}) + \\
\sum_{0\le|\b|<|\a|}
A_\b(Z(x),\1{Z(x)}) (L^\b\rho_Z)(Z(x),\1{Z(x)}) +
O(|x|^{k-l-1}).
\end{multline}
By the choices of $\rho(Z,\1Z)$ and $\rho'(Z',\1Z')$, we have the
decompositions in $\bC^{N_1}\times \bC^{N_2}$
\begin{equation}
L^\a\rho^j_Z=(L^\a\rho^j_{Z^1},0),\quad
L'{}^\a\rho'{}^j_{Z'}=(L'{}^\a\rho'{}^j_{Z'{}^1},0),
\quad j=1,\ldots,d.
\end{equation}
We multiply both sides of (\ref{kid3}) on the
right by the $N\times N_2$ constant matrix
$C=\left(0\atop{I}\right)$ with $I$ being the
$N_2\times N_2$ identity  matrix.
We conclude that
\begin{equation}\Label{kid4}
(L'{}^\a\rho'_{Z'{}^1})(Z'{}(x),\1{Z'{}(x)})
H^1_Z(Z(x)) C =  O(|x|^{k-l-1}).
\end{equation}
We now use the assumption that $M'_1$ is
$l$-nondegenerate.
By this assumption, we can choose multi-indices
$\a^1,\ldots,\a^{N_1}$ and integers
$j_1,\ldots,j_{N_1}$, with $0 \le |\a^\mu|\le l$,
$1 \le j_\mu\le d$, such that the
$N_1\times N_1$ matrix given by
$$B(x)
:=\Big(L'{}^{\a^\mu}\rho'{}^{j_\mu}_{Z'{}^1}
(Z'{}(x),\1{Z'{}(x)})\Big)
_{1\le \mu \le N_1}$$
is invertible for $x$ near $0$. Since
$B(x) H^1_Z(Z(x))C = O(|x|^{k-l-1})$ by
(\ref{kid4}) and
$H^1_{Z^2}(Z(x))\equiv H^1_Z(Z(x)) C$,
we conclude that $H^1_{Z^2}(Z(x))=
O(|x|^{k-l-1})$.  Since $M =M_1\times\C^{N_2}$ is
generic in
$\C^N=\C^{N_1}\times\C^{N_2}$,  the statement (i)
follows.

       From (i) it follows in particular that
$H^1_{Z^2}(0)=0$. Since $H$ is invertible, we
conclude that
$H^1_{Z^1}(0)$ is also invertible, and (ii)
follows from  (\ref{kids})
        by taking $x^2=y^2=0$. This completes the proof
of Lemma {\ref{add}}.
\end{proof}

We now give the proof of Theorem \ref{2}
assuming that Theorem
\ref{tech} has been proved.  As mentioned in the
beginning of this section, the proof of Theorem
\ref{tech} will be given in the remaining
sections.

\begin{proof}[Proof of Theorem~{\rm \ref{2}}] Set $V:=V_1\cup V_2\cup
V_3\subset M$, where $V_1,V_2,V_3$ are defined by
  (\ref{V1}), (\ref{V2}) and
(\ref{V3}) respectively. Let $p\in
M\setminus V$, $M'$ a real-analytic submanifold
of $\bC^N$ and
$p'\in M'$. We may assume that $M$ and $M'$ have
the same dimension, since otherwise there is
nothing to prove. Let
$\k>1$ be fixed.
       If, for
some integer $s > 1$,
$(M,p)$ and
$(M',p')$ are not
$s$-equivalent, then we can take  $k=s$
to satisfy the conclusion of
Theorem~\ref{2}.

Assume for the rest of the proof that $(M,p)$ and
$(M',p')$ are $k$-equivalent for all
$k>1$. Without loss of generality we may
assume $p = p' = 0$.  We shall make use of
Proposition~\ref{const}.  Since
$M$ is CR, of minimum degeneracy, and of
minimum orbit codimension at $0$,
$M'$ is also CR, of minimum degeneracy, and of
minimum orbit codimension at $0$.  Furthermore, in the
notation of Proposition ~\ref{const}, we have
$r_j(0) = r'_j(0)$,
$j=1,2, 3$.  Hence we may apply
Proposition~\ref{dec} to both $(M,0)$ and
$(M',0)$ with the same integers $N_1,N_2,N_3$
to obtain the decompositions
\begin{equation}\Label{mdecomp}
M=M_1\times\C^{N_2}\times\{0\},\ \ \ \
M'=M'_1\times\C^{N_2}\times\{0\},
\end{equation}
      where both
decompositions are understood in the sense of
germs at $0$ in $\bC^N$. Since
$M'_1$ is finitely nondegenerate at $0$, there
exists an integer
$l\ge 0$ such that
$M'_1$ is $l$-nondegenerate at $0$.

Assume first that $M$ and $M'$ are generic at
$0$, i.e.\ $N_3=0$.    Then, for every
$k>l+1$, the conclusions of Lemma~\ref{add} hold.
By conclusion (ii) of that lemma, for every
$k$-equivalence $H=(H^1,H^2)$ between $(M,0)$ and
$(M',0)$, the map
\begin{equation}\Label{F} h\colon Z^1\mapsto
H^1(Z^1,0)
\end{equation} is a $k$-equivalence between
$(M_1,0)$ and
$(M'_1,0)$. Furthermore, $(M_1,0)$ and $(M'_1,0)$
satisfy the assumptions of Theorem~\ref{tech}.

By  Theorem~\ref{tech}, there exists a
biholomorphic equivalence $\4h$ between
$(M_1,0)$ and $(M'_1,0)$ with
$\4h(Z^1)= h(Z^1)+O(|Z^1|^\k)$.
As we mentioned in the beginning of \S
\ref{kequiv}, without loss of generality, we can
assume that
$H$ is convergent.
Then we may define the germ
$\widehat H\colon(\C^N,0)\to(\C^N,0)$ of a
biholomorphism at the origin as follows:
\begin{gather}
\widehat H^1(Z^1,Z^2) := \widehat h(Z^1)  \\
\widehat H^2(Z) := H^2(Z)
\end{gather}
It is then a consequence of
Lemma~\ref{add}  that
$\widehat H$ satisfies the conclusion of
Theorem~\ref{2} if $k > \k + l + 1$.

We now return to the general case in which $M$ and
$M'$ are not necessarily generic, and let
$H=(H^1,H^2,H^3)$ be a $k$-equivalence corresponding
to the decomposition given by (\ref {mdecomp})
with $k > \k+l+1$. By Lemma
\ref{add0}  the mapping $(Z^1,Z^2) \mapsto
\big(H^1(Z^1,Z^2,0), H^2 (Z^1,Z^2,0)\big)$ is a
$k$-equivalence between the generic submanifolds
$(M_1\times\C^{N_2},0)$ and
$(M'_1\times\C^{N_2},0)$ and $H^3(Z^1,Z^2,0) =
O(|Z^1,Z^2|^k)$.  It follows from the generic
case, treated above, that there exists a
biholomorphic equivalence
$\4h(Z^1,Z^2)$ between
$(M_1\times\C^{N_2},0)$ and
$(M'_1\times\C^{N_2},0)$ such that
$\4h
(Z^1,Z^2) =  \big(H^1(Z^1,Z^2,0), H^2
(Z^1,Z^2,0)\big)+ O(|Z^1,Z^2|^\k)$.
We write
$\4h=(\4h^1, \4h^2)$ corresponding to the
product $\bC^{N_1}\times \bC^{N_2}$.  We may now
define
$\widehat H(Z^1,Z^2,Z^3)$ by
\begin{gather}
\widehat H^1(Z^1,Z^2,Z^3) := \widehat
h^1(Z^1,Z^2)+H^1(Z) - H^1(Z^1,Z^2,0)
\\
\widehat H^2(Z) :=
H^2(Z)\\
\widehat H^3(Z) := H^3(Z) -H^3(Z^1,Z^2,0)
\end{gather} and conclude that $\widehat H$
satisfies the desired conclusion of the theorem.
This completes the proof of Theorem \ref{2}
(assuming Theorem \ref{tech}).
\end{proof}

\section{Rings $\cR(V,V_0)$ of germs of holomorphic
functions}\Label{rings}

An important idea of the proof of Theorem~\ref{tech} is to
parametrize all $k$-equivalences
between $(M,0)$ and $(M',0)$ by their jets in an
expression of the form (\ref{gam}) below. For this
we shall introduce some notation for certain rings of
germs of holomorphic functions. 
If $V$ is a finite dimensional complex vector space
and $V_0\subset V$ is a vector subspace, we define
$\cR(V,V_0)$ to be the ring of all germs of
holomorphic functions
$f$ at $V_0$ in $V$ such that the restrictions
$\d^\alpha f|_{V_0}$ of all partial derivatives  are
polynomial functions on $V_0$. Here $\d^\a$ denotes
the partial derivative with respect to the
multiindex $\a\in\bZ_+^{\dim V}$ and some linear coordinates $x\in V$.
Recall that, if $f$ and
$g$ are two functions holomorphic in some
neighborhoods of
$V_0$ in $V$, then $f$ and $g$  define the same germ
of a holomorphic function at $V_0$ in $V$ if they
coincide in some (possibly smaller) neighborhood of
$V_0$ in $V$. We shall identify such a germ with any
representative of it. It is easy to see that the ring
$\cR(V,V_0)$ does not depend on the choice of linear coordinates 
in $V$ and is invariant under partial
differentiation with respect to these coordinates.
In the following we fix a complement $V_1$ of $V_0$
in $V$ so that we have the identification $V\cong
V_0\times V_1$ and fix linear coordinates
$x=(x_0,x_1)\in V_0\times V_1$.
In terms of these
coordinates  we may write an element $f\in
\cR(V,V_0)$ in the form
\begin{equation}\Label{ps} f(x_0,x_1)=\sum_\b
p_\b(x_0) x_1^\b, \quad
\b\in\bZ_+^{\dim V_1},
\end{equation} where the $p_\b(x_0)$ are polynomials
in $V_0$ satisfying the estimates
\begin{equation}\Label{pa} |p_\b(x_0)|\le
C(x_0)^{|\b|+1} \quad \text{for all } \b,
\end{equation} where $C(x_0)$ is a  positive locally
bounded function on
$V_0$. Conversely, every power series of the form
(\ref{ps}) satisfying (\ref{pa}) defines a unique
element of
$\cR(V,V_0)$.

In the following we shall consider germs of holomorphic maps
whose components are in $\cR(V,V_0)$. If $W$ is
another finite dimensional complex vector space and $W_0\subset W$ is a subspace,
     we shall write
$\phi\colon (V,V_0)\to (W,W_0)$ to mean
a germ at $V_0$ of a holomorphic map from $V$ to $W$ such that
$\phi(V_0)\subset W_0$.
It can be shown using the chain rule
that a composition $f\circ\phi$ with $\phi$ as 
above with components in $\cR(V,V_0)$
and $f\in \cR(W,W_0)$ always belongs to $\cR(V,V_0)$.
We shall prove the analogue of this property for
more general expressions which we shall need in the
 proof of Theorem \ref{iterth} below.

\begin{Lem}\Label{compose}  Let $V_0$, $V_1$,
$\2V_0$, $\2V_1$ be finite dimensional complex
vector spaces with fixed bases and $x_0$,
$x_1$, $\2x_0$, $\2x_1$ be the linear coordinates with
respect to these bases.
        Let
$q\in\C[x_0]$ and
$\2q\in\C[\2x_0]$ be nontrivial polynomial functions
on $V_0$ and
$\2V_0$ respectively, and let
$$\phi=(\phi_0,\phi_1)\colon \big(\C\times V_0\times
V_1,\C\times V_0\big)\to
\big(\2V_0\times\2V_1,\2V_0\big)$$ be a germ of a
holomorphic map with components in the ring
$\cR\big(\C\times V_0\times V_1, \C\times V_0\big)$
and satisfying
\begin{equation}\Label{nonvanishing}
\2q\Big(\phi_0\big(\frac{1}{q(x_0)},x_0,0\big)\Big)
\not\equiv 0.
\end{equation}
Then there exists a
ring homomorphism
\begin{equation}\Label{hom}
\cR\big(\C\times\2V_0\times\2V_1,\C\times\2V_0\big)\ni
\2f\mapsto f\in \cR\big(\C\times V_0\times V_1,
\C\times V_0\big)
\end{equation}
such that
\begin{equation}\Label{fft}
\2f\bigg(\frac{1}
{\2q\big(\phi_0(\frac{1}{q(x_0)},x_0,x_1)\big)},
\phi\Big(\frac{1}{q(x_0)},x_0,x_1\Big)\bigg) \equiv
f\Big(\frac{1}{p(x_0)},x_0,x_1\Big),
\end{equation} 
 with
$p(x_0):=q(x_0)^{d_0+1}
\2q\big(\phi_0(\frac{1}{q(x_0)},x_0,0)\big)$,
where $d_0$ is the degree of the polynomial
$(\theta,x_0)\mapsto
\2q\big(\phi_0(\theta,x_0,0)\big)$ with
respect to $\theta$.  
Furthermore,  $f$ vanishes on $\C\times V_0$ if $\2f$ vanishes
on $\C\times\2V_0$.
\end{Lem}

\begin{proof} For $\2f$ as above and $\T',\T''\in\C$, define a germ
$g$ at $\C\times\C\times V_0$
of a holomorphic function on $\C\times\C\times V_0\times V_1$ by
\begin{equation} g(\T',\T'',x_0,x_1):=
\2f\Big(\frac{\T'} {1+\T'[\2q\big(\phi_0(\T'',x_0,x_1)\big)
-
\2q(\phi_0(\T'',x_0,0))]},
\phi(\T'',x_0,x_1)
\Big).
\end{equation} 
We use the consequence of the chain rule
that any partial derivative of a composition
of two holomorphic maps can be written as a polynomial expression
in the partial derivatives of the components.
Then it follows from the assumptions of
the lemma that $g$ is in the ring
$\cR\big(\C\times\C\times V_0\times
V_1,\C\times\C\times V_0\big)$. It is straightforward
to see that, if $f$ is given by
$$f(\T,x_0,x_1):=g\Big(\T q(x_0)^{d_0+1}, \T
q(x_0)^{d_0}
\2q\big(\phi_0(\frac{1}{q(x_0)},x_0,0)\big),x_0,x_1
\Big),$$
then (\ref{fft}) holds and the map
$\2f\mapsto f$ satisfies the conclusion of the lemma.
\end{proof}

\section{Jet spaces of mappings}
\Label{equivs}
For integers $r,m,l\ge 0$, we denote by
$J^r_{m,l}$ the space of all jets at $0$ of order $r$ of
holomorphic maps from $\C^m$ to $\C^l$. This
is a complex vector space that can be identified
with the space of
$\C^l$-valued polynomials on $\C^m$ of degree at
most $r$.  We write such a polynomial in the form
$\sum_{0\le|\a|\le r}(\lambda_\a/\a!) Z^\a$,
$\lambda_\a\in\C^l$, and call
$(\lambda_\a)_{0\le |\a|\le r}$ the standard
linear coordinates in $J^r_{m,l}$. For fixed integers
$n,d\ge 0$ and $N:=n+d$, we introduce the complex
vector spaces
\begin{equation}\Label{er} E^{r}:=J^r_{N,N}\times
J^r_{n,d}\times\bC^n,
\quad E^{r}_0:=J^r_{N,N}\times\{(0,0)\},
\quad E^{r}_1:=\{0\}\times J^r_{n,d}\times\bC^n
\end{equation}
with $E^{r}_0,E^{r}_1\subset E^r$.  We use the
identification
$E^r\cong E^{r}_0\times E^{r}_1$.

Let $M$ and $M'$ satisfy the assumptions of
Theorem~\ref{tech}. According to
Proposition~\ref{dec-orbit} we write $M$ near $0$ in the form
\begin{equation}\Label {mq}
M=\big\{(z,w,u)\in\bC^n\times\bC^{d_1}
\times\bR^{d_2}:
w=Q(z,\1z,\1w,u)\big\},
\end{equation} where $Q$ is a germ at $0$
in
$\bC^{2n+d}$ of a holomorphic
$\bC^{d_1}$-valued function satisfying conditions
(\ref{Qcond}). We also choose normal coordinates for
$M'$ so that
$$M'=\big\{(z',w')\in\bC^n\times\bC^d:
w'=Q'(z',\1{z'},\1{w'})\big\},$$ where $Q'$ is a germ
at
$0$ in
$\bC^{2n+d}$ of a holomorphic $\bC^d$-valued function
satisfying
\begin{equation}\Label{Q'cond} Q'(z',0,\w')\equiv
Q'(0,\z',\w')\equiv\w'.
\end{equation}

In these coordinates (which will be fixed for the
remainder of the paper), for
every invertible germ of a holomorphic map
$H\colon (\bC^N,0)\to (\bC^N,0)$ we write
$H(Z)=(F(Z),G(Z))$ with
$z'=F(Z)$,
$w'=G(Z)$ and
$Z=(z,w,u)$.
For $Z\in\C^{N}$ near the origin, we define
\begin{equation}\Label{jr}
\D^rH(Z):=
\bigg(
\, \Big(\frac{\d^{|\a|} H}{\d
Z^\a}(Z)\Big)_{0\le|\a|\le r},
\, \Big(\frac{\d^{|\nu|} G}{\d
z^\nu}(Z)\Big)_{0\le|\nu|\le r},
\,F(Z)\bigg).
\end{equation} We think of $\D^rH$ as a germ at $0$
of a holomorphic map from $\bC^{N}$ into  the vector
space $E^r$ defined by (\ref{er}).

Now let $H=(F,G)$ be a $k$-equivalence between
$(M,0)$ and $(M',0)$ with $k > 1$.
By a standard complexification argument, $H$ is a
$k$-equivalence  means that the identity
\begin{equation}\Label{qid}  G(z,Q(z,\z,\w,u),u)
\equiv
Q'\left(F(z,Q(z,\z,\w,u),u),\Hb(\z,\w,u)\right) +
R(z,\z,\w,u)
\end{equation} holds for all
$(z,\z,\w,u)\in\bC^{2n+d}$ near the origin, where
$R(z,\z,\w,u) = O(k)$. In particular, for
$(\z,\w,u)=0$ we obtain from (\ref{qid}) and
(\ref{Q'cond}) the identity
\begin{equation}\Label{Gid}
G(z,0)= O(k)
\end{equation}
and hence for $r < k$, we have
$\D^rH(0)\in E_0^r$, where $E_0^r$ is defined by
(\ref{er}).

\section{The basic identity} We assume that the
assumptions of Theorem \ref{tech} hold and that
$\rho(Z,\1Z)$ is a defining function for $M$ at $0$.
We begin by establishing a relation, called the
{\it basic identity}, between two jets of a
$k$-equivalence
$H$ at  points
$Z$ and $\Z$ in $\C^N$ satisfying $(Z,\Z)\in\3M$,
i.e. $\rho(Z,\Z)=0$. We shall make use
of the notation introduced in
\S~\ref{rings}-\ref{equivs}. In particular, we
have normal coordinates
$(z,w,u)\in\C^n\times\C^{d_1}\times\C^{d_2}$ for $M$
and
$(z',w')\in\C^n\times\C^d$ for $M'$ and write
$Z=(z,w,u)$, $\Z=(\z,\w,u)$. Furthermore we use
matrix notation and regard
$F_z(Z)$ as an $n\times n$ matrix,
$F_w(Z)$ as an ${n\times d_1}$ matrix,
$G_z(Z)$ as a ${d\times n}$  matrix  and $G_w(Z)$ as
a ${d\times d_1}$ matrix.  Similarly
$Q_z(z,\z,\w,u)$ is regarded as a $d_1\times n$ matrix.

To shorten the notation it will be convenient to
write for $r,m$ nonnegative integers
\begin{equation}\Label{rm}
\cR^{r}_{m}:=\cR\big(\C\times E^{r}\times \C^{m},
\C\times E^{r}_0\times \{0\}\big),
\end{equation}
	 where the rings $\cR(V,V_0)$ are defined as in \S
\ref{rings} and the vector spaces $E^{r}$ and
$E^{r}_0$ are defined in (\ref{er}). We can now
state precisely the basic identity.

\begin{Thm}[Basic Identity]\Label{basic} Let
$(M,0)$ and $(M',0)$ be two germs of generic
real-analytic submanifolds of $\bC^N$ satisfying the
assumptions of Theorem \ref{tech}.  Assume that
$M'$ is
$l$-nondegenerate at $0$ (with $l\ge 0$) and that
normal coordinates for $M$ and $M'$ are chosen as
above. Then for every integer
$r>0$, there exists a germ of a holomorphic map
\begin{equation}\Label{van}
\Psi^r\colon \big(\C\times E^{r+l}\times \C^{2N},
\C\times E^{r+l}_0\times \{0\}\big)\to (E^r,E^r_0)
\end{equation}
and for $r=0$, a germ
$\Psi^0\colon \big(\C\times E^l\times \C^{2N},
\C\times E^l_0\times \{0\}\big)\to (E^0,0)$,
such that the components of $\Psi^r$, $r\ge 0$, are in the ring
$\cR^{r+l}_{2N}$ and the following holds.
For every
$k$-equivalence $H=(F,G)$ between $(M,0)$ and
$(M',0)$ with
$k>r+l$, one has for $(Z,\Z)$ near the origin in
$\bC^{2N}$,
\begin{equation}\Label{desired}
\D^r H(Z) = \Psi^r\bigg(\frac{1}{\det(\1F_\z(\Z))},
\D^{r+l}\1H(\Z),\Z,Z\bigg) + R_H^r(Z,\Z),
\end{equation}
where $R_H^r(Z,\Z)$  is a germ at
$0$ of a holomorphic map from $\C^{2N}$ into $E^r$
whose restriction to
$\3M$ vanishes of order $k-r-l$ at $0$.
\end{Thm}

\begin{proof}  For convenience we use the notation
$$\omega:=(z,\Z)=(z,\z,\w,u)\in
\C^n\times\C^n\times\C^{d_1}\times\C^{d_2},\quad
Z(\o):=(z,Q(\omega),u) \in \C^n\times\C^{d_1}\times\C^{d_2},$$
so that the equation of
$\3M\subset\C^{2N}$ near $0$ is given by $w
= Q(\o)$, or equivalently, by $Z=Z(\o)$. We first
differentiate the identity (\ref{qid}) in
$z\in\bC^n$. Using the chain rule we obtain the
identity in matrix notation
\begin{multline}\Label{id1}
G_z(Z(\o)) + G_w(Z(\o))
Q_z(\o) \equiv \\ Q'_{z'}\left(F(Z(\o)),\Hb(\Z)\right)
\big(F_z(Z(\o)) + F_w(Z(\o)) Q_z(\o)\big) + R_z(\o),
\end{multline} where $R_z(\o)=O(|\o|^{k-1})$. (Observe that
$R_z$ in (\ref{id1}) depends on the map $H$). The
invertibility of
$H$ implies the invertibility of $F_z(0)$ and hence
of
$F_z(Z(\o)) + F_w(Z(\o)) Q_z(\o)$ for $\o$ near the
origin (since
$Q_z(0)=0$ by (\ref{Qcond})). Hence we conclude for
$\o$ sufficiently small,
\begin{multline}\Label{qzid} Q'_{z'}
\left(F(Z(\o)),\Hb(\Z)\right) = \\
\big(G_z(Z(\o)) + G_w(Z(\o))Q_z(\o)\big)
\big(F_z (Z(\o))+ F_w(Z(\o)) Q_z(\o)\big)^{-1} +
O(|\o|^{k-1}).
\end{multline} 
Our next goal will be to express the
right-hand side of (\ref{qzid}) and then its derivatives
in terms of functions in $\cR^{r}_{2n+d}$
that vanish on certain vector subspaces.
For this we introduce the notation
\begin{multline}\Label{br} A^r:=\bC\times
(J^r_{N,N}\times\{0\}\times
\bC^n)\times(\bC^n\times\{0\}) \subset  \\
\bC\times (J^r_{N,N}\times J^r_{n,d}\times
\bC^n)\times(\bC^n\times\C^{n+d}) = \C\times
E^r\times\C^{2n+d}.
\end{multline}
We have the following lemma.
\begin{Lem}
\Label{GF} With the notation above there exists  a
$d\times n$ matrix $P$, independent of $H$,
        with entries in
$\cR^{1}_{2n+d}$ such that, for  $\o$
in a neighborhood of $0$ in $\C^{2n+d}$,
\begin{multline}\Label{pqz}
\big(G_z(Z(\o)) + G_w(Z(\o)) Q_z(\o)\big)
\big(F_z(Z(\o))+ F_w(Z(\o))Q_z(\o)\big)^{-1}\equiv\\
P\Big(\frac{1}{\det
F_z(Z(\o))},\D^1H(Z(\o)),\o\Big)
\end{multline}  and $P$ vanishes on the subspace
$A^1\subset\C\times E^1\times\C^{2n+d}$ defined by 
{\rm(\ref{br})}.
\end{Lem}

\begin{proof} For simplicity we drop the argument
$Z(\o)$ in $G_z$,
$G_w$, $F_z$, $F_w$ and $\D^1H$. We have
\begin{equation}\Label{pqz11}
\big(G_z + G_w  Q_z(\o)\big)
\big(F_z + F_w  Q_z(\o)\big)^{-1}\equiv
\big(G_z  + G_w  Q_z(\o)\big)
\big(I+ F_z ^{-1}F_w Q_z(\o)\big)^{-1}F_z^{-1}.
\end{equation}
The first factor in the right-hand
side of (\ref{pqz11}) can be expressed as a matrix
valued polynomial in the entries of $G_z$ and $G_w$
with holomorphic coefficients in $\o$. We now think
of the entries of $G_z$ as variables in $J^1_{n,d}$
and those of $G_w$ as part of the variables in
$J^1_{N,N}$ and write
$$\big(G_z  + G_w
Q_z(\o)\big)\equiv P_1\Big(\frac{1}{\det
F_z},\D^1H,\o\Big)$$
with $P_1$  independent of the variable in the
first factor $\C$ and having entries
in $\cR^{1}_{2n+d}$.  Since $Q_z(z,0,0,0)\equiv 0$,
$P_1$ vanishes on the subspace
$A^1\subset\C\times E^1\times\C^{2n+d}$ defined by
(\ref{br}) with $r=1$. By the standard formula for
the inverse of a matrix, the third factor in the
right-hand side of (\ref{pqz11})
       can be also written in the form
$P_3\big(\frac{1}{\det F_z},\D^1H,\o\big)$, where
$P_3$ is a matrix valued polynomial (with
entries in $\cR^{1}_{2n+d}$) depending only on part
of the variables in $\C\times J^1_{N,N}$ and
independent of the variables in
$J^1_{n,d}\times\C^n$ and $\o$. The second
factor in the right-hand side of (\ref{pqz11}) can
also be written in the form
$P_2\big(\frac{1}{\det F_z},\D^1H,\o\big)$ with
the entries of $P_2$ in $\cR^{1}_{2n+d}$. This can
be shown by using the chain rule in addition to the
arguments used for the first and third factors.  The
proof of the lemma is completed by taking
$P:=P_1 P_2 P_3$ and using the fact that
$\cR^{1}_{2n+d}$ is a ring.
\end{proof}

For the sequel we shall need the following lemma, which
is  proved by  repeated use of the chain rule, making
use of the identities (\ref{Gid}),
$Q(z,0,0,0)\equiv 0$, and induction on
$|\a|$. The details are left to the reader.

\begin{Lem}\Label{jump} Let $M$ and $M'$ be as
in Theorem {\rm\ref{basic}}. Then for every $f\in
\cR^{r}_{2n+d}$ with
$r\ge 1$ and every
$\a\in\bZ_{+}^{2n+d}$, there exists
$f^\alpha\in\cR^{r+|\a|}_{2n+d}$ such that the
following holds.  For any $k$-equivalence
$H=(F,G)$ between
$(M,0)$ and $(M',0)$  with $k > r+|\a|$,
\begin{multline}\Label{dif}
\big({\d^{|\a|}}/{\d \o^\a}\big)
f\Big(\frac{1}{\det
F_z(Z(\o))},\D^{r}H(Z(\o)),\o\Big) \equiv\\
f^\alpha\Big(\frac{1}{\det
F_z(Z(\o))},\D^{r+|\a|}H(Z(\o)),\o\Big).
\end{multline}   If in addition
$\a\in\bZ_{+}^n\times \{0\}$ (i.e. the
differentiation in {\rm(\ref{dif})} is taken with respect to
$z$ only) and if
$f$ vanishes on the subspace
$A^r\subset \C\times E^r\times\C^{2n+d}$
defined by
{\rm (\ref{br})}, then $f^\a$ vanishes on the subspace
$A^{r+|\a|}\subset \C\times
E^{r+|\a|}\times\C^{2n+d}$.
\end{Lem}

We now return to the proof of Theorem~\ref{basic}.
By making use of (\ref{qzid}) and (\ref{pqz}) we
obtain the identity
\begin{equation}\Label{pqz1}
Q'_{z'}\left(F(Z(\o)),\Hb(\Z)\right) =
P\Big(\frac{1}{\det
F_z(Z(\o))},\D^1H(Z(\o)),\o\Big) + O(|\o|^{k-1}),
\end{equation} where $\Z=(\z,\w,u)$ as before and $P$ is given
by Lemma~\ref{GF}.

We claim that for every
$\b\in\bZ_{+}^n$ with $0\le |\b|\le l$, there exists
$P^\b\in (\cR^{|\b|}_{2n+d})^d$, independent of
$H$, vanishing on the subspace
$A^{|\b|} \subset  \bC\times
E^{|\b|}\times\C^{2n+d}$, and such that the following
identity holds for $\o$ in a neighborhood of $0$ in
$\C^{2n+d}$:
\begin{equation}\Label{qa}
Q'_{z'{}^{\b}}\left(F(Z(\o)),\Hb(\Z)\right) =
P^\b\Big(\frac{1}{\det
F_z(Z(\o))},\D^{|\b|}H(Z(\o)),\o\Big) +
O(|\o|^{k-|\b|}).
\end{equation}
Indeed, for $\b=0$, (\ref{qa})
follows directly from  (\ref{qid}) and for $|\b|=1$,
(\ref{qa}) is a reformulation of (\ref{pqz1}). For
$|\b|>1$ we prove the claim by induction on $|\b|$.
Assume that (\ref{qa}) holds for some $\b$.  By
differentiating (\ref{qa}) with respect to $z$ we
obtain in matrix notation the identity
\begin{multline}\Label{qzid1}
\big(Q'_{z'{}^{\b}}\big)_{z'}
\left(F(Z(\o)),\Hb(\Z)\right) \big(F_z (Z(\o))+
F_w(Z(\o)) Q_z(\o)\big)=\\
\left(\d/\d z\right) P^\b\Big(\frac{1}{\det
F_z(Z(\o))},\D^{|\b|}H(Z(\o)),\o\Big) + O(|\o|^{k-|\b|-1}).
\end{multline}
By Lemma~\ref{jump} we have
\begin{equation}\Label{der}
\left(\d/\d z\right) P^\b\Big(\frac{1}{\det
F_z(Z(\o))},\D^{|\b|}H(Z(\o)),\o\Big) \equiv
S\Big(\frac{1}{\det
F_z(Z(\o))},\D^{|\b|+1}H(Z(\o)),\o\Big),
\end{equation}
where $S$ is a $d\times n$ matrix
with entries in
$\cR^{|\b|+1}_{2n+d}$, vanishing on the subspace
$A^{|\b|+1}\subset\C\times E^{|\b|+1}\times\C^{2n+d}$.
Since, as in  the proof of Lemma~\ref{GF}, each
entry of the matrix $\big(F_z (Z(\o))+ F_w(Z(\o))
Q_z(\o)\big)^{-1}$ can be written in the form
$f\big(\frac{1}{\det F_z(Z(\o))},\D^{1}H(Z(\o)),\o\big)$
with $f$ in
the ring
$\cR^{1}_{2n+d}$, the identity (\ref{qa})  for
$\b$ replaced by any multiindex $\b'$ with
$|\b'|=|\b|+1$ follows
from (\ref{qzid1}) and (\ref{der}) by observing that the ring
$\cR^{1}_{2n+d}$ has a natural embedding into
$\cR^{|\b|+1}_{2n+d}$. This completes the proof of
the claim.

We now use the condition that $M'$ is
$l$-nondegenerate which is equivalent to
$$\span_\C \big\{  Q'{}^j_{z'{}^\b\z'}(0,0,0) : 1\le j\le
d,\, 1\le |\b| \le l \big\} = \bC^n$$
(see e.g. \cite{BER}, 11.2.14).
    From this,
together with (\ref{Q'cond}), we conclude that  we
can select a subsystem of $N$ scalar identities from
(\ref{qa}) from which
$\Hb(\Z)$ can be solved uniquely by the implicit
function theorem. We obtain
\begin{equation}\Label{H}
\Hb(\Z) = T\bigg(F(Z(\o)),
P^{\b}\Big(\frac{1}{\det
F_z(Z(\o))},\D^{|\b|}H(Z(\o)),\o\Big)_{0\le|\b|\le
l}
\bigg)  +
  O(|\o|^{k-l}),
\end{equation}
       where $T$ is a germ of a holomorphic
map $T\colon (\C^n\times\C^m,0)\to (\C^N,0)$, with
\begin{equation}\Label{m} m:= d\times \#
\{\b\in\bZ^n_+ : 0\le|\b|\le l\}.
\end{equation}
        Observe that the germ $T$ depends only on
$Q'$ but not on $H$.

We claim that there exists
$\Phi\in (\cR^{l}_{2n+d})^N$, independent of $H$, such that
\begin{equation}\Label{Hclaim}
\Hb(\Z) = \Phi\Big(\frac{1}{\det
F_z(Z(\o))},\D^lH(Z(\o)),\o\Big) + O(|\o|^{k-l}).
\end{equation}

In order to prove the claim we use the notation
$x_0 := (\T,\L)\in \bC\times J^l_{N,N}$ and $x_1:=
(\L',z',\o) \in J^l_{n,d}\times\bC^n\times \bC^{2n+d}$,
and for $l \ge r$, we denote by $\pi_r^l:E^l\to E^r$
the natural projection from  $E^l$ onto $E^r$. We
define $\Phi$ by
\begin {equation}\Label{PhiT}
\Phi(\T,\L,\L',z', \o): = T\bigg(z',
P^\b\big (\T,
\pi^l_{|\b|}(\L,\L',z'),\o\big)_{0\le |\b|\le
l}\bigg).
\end{equation}
To show that $\Phi$ is in $(\cR^{l}_{2n+d})^N$, we must
differentiate the right hand side of (\ref{PhiT}) with
respect to $x_1=(\L',z',\o)$ and evaluate at $x_1=0$.
By using the chain rule and the fact that  each $P^\b$
is in
$(\cR^{|\b|}_{2n+d})^d$ and vanishes when $x_1 = 0$,
it is easy to check that for any multiindex $\a$,
$${\partial^\a\over \partial x_1^\a}T\bigg(z',P^\b\big
(\T,
\pi^l_{|\b|}(\L,\L',z'),\o\big)_{0\le |\b|\le
l}\bigg)\bigg|_{x_1 = 0}$$ is a polynomial in $x_0$.
This proves the claim (\ref{Hclaim}).

We now differentiate the identity (\ref{Hclaim})
with respect to
$\Z=(\z,\w,u)$. By using Lemma~\ref{jump} again, we find
$\Phi^\b\in(\cR^{l+|\b|}_{2n+d})^N$, independent of
$H$, such that
\begin{equation}\Label{desired1}
\d^\b\ov H(\Z) = \Phi^\b\Big(\frac{1}{\det
F_z(Z(\o))},\D^{l+|\b|}H(Z(\o)),\o\Big) +
O(|\o|^{k-l-|\b|}).
\end{equation} For any $\b\in\bZ_{+}^N$ we decompose
$\Phi^\b = (\Phi_1^\b,\Phi_2^\b)\in \bC^n\times\bC^d$
and set for $\T, \L,\L', z',\o$ as above and
$Z=(z,w,v)\in\C^n\times\C^{d_1}\times\C^{d_2}$,
\begin{multline}\Label{}
\2\Phi^\b(\T,\L,\L',z',Z,\Z):=\\
\begin{cases}
\Phi^0(\T,\L,\L',z',\o)-\Phi^0(\T,\L,0,0,0) &
\text{ for } \b=0,\\
\big(\Phi_1^\b(\T,\L,\L',z',\o),\Phi_2^\b(\T,\L,\L',z'
,\o)-
\Phi_2^\b(\T,\L,0,0,0)\big)  & \text{ for }
\b \in\bZ_{+}^n\times\{0\},\b\neq 0\\
\Phi^\b(\T,\L,\L',z',\o) &\text{ otherwise }.
\end{cases}
\end{multline}
Clearly
$\2\Phi^\b$ is in $(\cR^{l+|\b|}_{2N})^N$ and is
independent of
$w$ and $v$.
Since for any $k$-equivalence $H=(F,G)$ we have
$\d^\b G(0) = 0$ for $\b\in\bZ_{+}^n\times\{0\}$ with
$|\b|< k$ by (\ref{Gid}), it follows from (\ref{desired1}) that
$\2\Phi_2^\b\big(\frac{1}{\det
F_z(0)},\D^{l+|\b|}H(0),0\big)= 0$.
Hence (\ref{desired1}) implies
\begin{equation}\Label{desired2}
\d^\b\ov H(\Z) = \2\Phi^\b\left(\frac{1}{\det
F_z(Z)},\D^{l+|\b|}H(Z),Z,\Z\right) +
\2 R_H^\b(Z,\Z),
\end{equation} where $R_H^\b$ is a germ  at $0$ of a
       holomorphic map from $\C^{2N}$ to $\C^N$ depending on $H$ and
whose restriction to $\3M$ vanishes  at $0$ of order
$k-l-|\b|$. By taking complex conjugates of
(\ref{desired2}) for $0\le|\b|\le r$, and using the fact
that $(Z,\Z)\in\3M$ is equivalent to
$(\1\Z,\1Z)\in\3M$, we obtain (\ref{desired}) with
$\Psi^r$ satisfying the conclusion of Theorem~\ref{basic}.
\end{proof}


\section{The iterated basic identity}
In this section we apply the relation given by
Theorem~\ref{basic} to different points and
iterate them, i.e. substitute one into the next
and so on. Let $(M,0)$ and $(M',0)$ satisfy the
assumptions of Theorem~\ref{basic}.  If
$\rho(Z,\1Z)$ is a defining function of $M$ near
$0$ and $s\ge 1$ is an integer, we define  a germ
$\3M^{2s}$ at $0$ of a complex manifold of
$\C^{(2s+1)N}$ by
\begin{multline}\Label{m2s}
\3M^{2\s}:=\big\{ (Z,\Z^1,Z^1,\ldots,\Z^\s,Z^\s)\in
\C^{(2\s+1)N} : \\
\rho(Z,\Z^1)=\cdots=\rho(Z^{\s-1},\Z^\s)=
\rho(Z^1,\Z^1)=\cdots=\rho(Z^s,\Z^s)=0
\big\}.
\end{multline} Hence $\3M^{2s}$ has codimension
$2sd$ in $\C^{(2s+1)N}$, where $d$ is the codimension
of $M$ in $\C^N$. (The iterated complexification
$\3M^{2s}$ was introduced by the third author in
\cite{Zait}.  For $Z_s$ fixed in (\ref{m2s}), this
corresponds to the Segre manifold of order $2s$ of $M$ at
$Z_s$ in the terminology of \cite{BER}.) For a
$k$-equivalence
$H$ between
$(M,0)$ and
$(M',0)$, we use the notation $\D^rH(Z)$ introduced
in (\ref{jr}). It will be also convenient to write
$$j^rH(Z):= \Big(\frac{\d^{|\a|} H}{\d
Z^\a}(Z)\Big)_{0\le|\a|\le r},$$ which is the first
$J^r_{N,N}$-valued component of $\D^rH(Z)$.  The
main result of this section is the following.

\begin{Thm}\Label{iterth}
        Under the assumptions of
Theorem~{\rm\ref{basic}}, for all integers $r\ge 0$ and
$s\ge 1$, there exists a polynomial $q^r_s$ on
$J^{r+2sl}_{N,N}$ and, for $r>0$, a germ
\begin{equation}\Label{}
\Psi^{r,s}\colon \big(\C\times E^{r+2sl}\times
\C^{(2s+1)N},
\C\times E^{r+2sl}_0\times\{0\}\big)\mapsto
(E^r,E^r_0)
\end{equation}
and for $r=0$, a germ $\Psi^{0,s}\colon
\big(\C\times E^{2sl}\times\C^{(2s+1)N}, \C\times
E^{2sl}_0\times\{0\}\big)\mapsto (E^0,0)$,
whose components are in the ring
$\cR^{r+2sl}_{(2s+1)N}$ such that, if $H=(F,G)$ is a
$k$-equivalence between $(M,0)$ and
$(M',0)$ with $k> 2sl+r$, the following holds:
\begin{equation}\Label{detprop}
q^r_s(j^{r+2sl}H(0))= \big(\det
F_z(0)\big)^{a^r_s}
\ov{\big(\det F_z(0)}\big)^{b^r_s},
\text{ for  some } a^r_s, b^r_s\in\bZ_+,
\end{equation}
\begin{multline}\Label{iter}
\D^r H(Z) =
\Psi^{r,s}\left(\frac{1}{q^r_s(j^{r+2sl}H(Z^s))},
\D^{r+2sl}H(Z^s),Z,\Z^1,Z^1,\ldots,\Z^s,Z^s\right)
+\\ R_H^{r,s}(Z,\Z^1,Z^1,\ldots,\Z^s,Z^s),
\end{multline}
where $R_H^{r,s}$ is a germ at $0$
of a holomorphic map
from $\C^{(2s+1)N}$ to
$E^r$, depending on $H$, whose restriction to
$\3M^{2s}$ vanishes of order
$k-r-2sl$ at $0$.
\end{Thm}
Note that since $H$ is a $k$-equivalence, it
follows that $\det F_z(0)\neq 0$, and hence the
right hand side of (\ref{detprop}) is necessarily
nonvanishing.
\begin{proof} We prove the theorem by induction on
$s\ge 1$. We start first with the case $s=1$ and
assume that $H=(F,G)$ is a $k$-equivalence between
$(M,0)$ and
$(M',0)$ with $k> r+2l$. By conjugating
(\ref{desired}) with
$r$ replaced by
$r+l$  we obtain
\begin{equation}\Label{desired12}
\D^{r+l} \1H(\Z) \equiv
\1{\Psi^{r+l}}\Big(\frac{1}{\det(F_z(Z^1))},
\D^{r+2l} H(Z^1),Z^1,\Z\Big) +
\1{R_H^{r+l}}(\Z,Z^1)
\end{equation}
with $\Psi^{r+l}$  and $R_H^{r+l}$
as in Theorem~\ref{basic}. If we observe that
$(Z,\Z)\in\3M\iff (\1\Z,\1Z)\in\3M$, we conclude
that the second term on the right hand side of
(\ref{desired12}) vanishes at $0$ of order $k-r-2l$
when
$(Z^1,\Z)\in\3M$.  Our next goal will be to
substitute (\ref{desired12}) into (\ref{desired})
and to apply Lemma~\ref{compose}. For this, we
define polynomials $q\in\C[\L]$ and
$\2q\in\C[\2\L]$, for
$\L=x_0\in V_0:=E^{r+2l}_0\cong J^{r+2l}_{N,N}$ and
$\2\L=\2x_0\in \2V_0:=E^{r+l}_0\cong J^{r+l}_{N,N}$, to be the
determinants of the parts of the jets $\L$ and
$\2\L$ obtained from the  first $n$ rows and first
$n$ columns of the linear terms of $\L$ and $\2\L$
respectively
(i.e. corresponding to $\det F_z(Z)$ and to $\det \1F_\z(\Z)$
for $\L=j^{r+2l}H(Z)$ and $\2\L=j^{r+l}\1H(\Z)$ respectively).
We also set
$x_1=(\L',z',Z,\Z,Z^1)\in V_1:=
J^{r+2l}_{n,d}\times\C^n\times\C^N\times\C^N\times\C^{N}
\cong E^{r+2l}_1\times\C^{3N}$,
$\2V_1:=E^{r+l}_1\times\C^{2N}$ and for $\T\in\C$,
$$\phi(\T,\L,\L',z',Z,\Z,Z^1):=\big(\1{\Psi^{r+l}}(\T,\L,\L',z',Z^1,\Z
),\Z,Z\big)\in
  E^{r+l}\times\C^{N}\times \bC^N.$$
(Observe that $E^{r+l}\times\C^{2N}=\2V_0\times \2V_1 $ by the
definition of $\2V_0$ and $\2V_1$ above.) Then $\phi$
satisfies the assumptions of Lemma~\ref{compose},
in particular, (\ref{nonvanishing}) holds since by
(\ref{desired12}) we have
\begin{equation}\Label{qtil}
\2q\bigg(
\1{\Psi_0^{r+l}}
\Big(\frac{1}{q(j^{r+2l}H(0))},
\D^{r+2l} H(0),0,0\Big)
\bigg) =\det \1F_\z(0),
\end{equation}
and the right hand side of (\ref{qtil}) is nonvanishing 
whenever $H=(F,G)$ is a $k$-equivalence with $k>1$.

    From substituting (\ref{desired12}) into
(\ref{desired}) we obtain the identity
\begin{multline}\Label{}
\D^r H(Z) \equiv
\Psi^r\bigg(
\frac{1} {\2q\big(
\1{\Psi_0^{r+l}}
(\frac{1}{q(j^{r+2l}H(Z^1))},
\D^{r+2l} H(Z^1),Z^1,\Z)
\big)},\\
\1{\Psi^{r+l}}
\Big(\frac{1}{q(j^{r+2l}H(Z^1))},
\D^{r+2l} H(Z^1),Z^1,\Z\Big),\Z,Z
\bigg) + R_H^{r,1}(Z,\Z,Z^1),
\end{multline}
where the restriction of $R_H^{r,1}$
to $\3M^2\subset\C^{3N}$ vanishes of order
$k-r-2l$ at the origin. Then for
$s=1$,  (\ref{iter}) is a consequence of Lemma~\ref{compose} with
$q^r_1$ being the polynomial $p$ given by the lemma.
The required property (\ref{detprop}) follows from
(\ref{qtil}) and from the explicit formula for $p$ in the
lemma.

Now we assume that (\ref{detprop}) and (\ref{iter}) hold for some fixed $s\ge
1$ and any $r\ge 0$ and shall prove them for $s+1$
and any $r\ge 0$. We replace the terms $j^{r+2sl}H(Z^s)$
and $\D^{r+2sl}H(Z^s)$
by using (\ref{iter}) with
$s=1$ and $r$ replaced by $r+2sl$.
We obtain
\begin{multline}\Label{}
\D^r H(Z) \equiv\\
\Psi^{r,s}\bigg(
\frac{1} {q^r_s\big(
\Psi^{r+2sl,1}_0
(\frac{1}{q^{r+2sl}_1(j^{r+2l(s+1)}H(Z^{s+1}))},
\D^{r+2l(s+1)} H(Z^{s+1}),Z^s,\Z^{s+1},Z^{s+1})
\big)},\\
\Psi^{r+2sl,1}
\Big(\frac{1}{q_1^{r+2sl}(j^{r+2l(s+1)}H(Z^{s+1}))},
\D^{r+2l(s+1)}
H(Z^{s+1}),Z^s,\Z^{s+1},Z^{s+1}\Big),\\
Z,\Z^1,Z^1,\ldots,\Z^s,Z^s
\bigg) +
R_H^{r,s+1}(Z,\Z^1,Z^1,\ldots,\Z^{s+1},Z^{s+1})
\end{multline}
with the restriction of $R_H^{r,s+1}$ to
$\3M^{2(s+1)}\subset\C^{(2(s+1)+1)N}$
vanishing of order $k-r-2l(s+1)$.
Similarly to the preceeding proof of (\ref{iter}) for $s=1$,
the desired conclusion of the theorem follows by
making use of Lemma~\ref{compose}.
\end{proof}

\section{Reducing the number of parameters}\label{gampar}
The expression in the right-hand side of (\ref{iter}) depends
on $(2s+1)N$ complex variables. Our goal in this section
will be to reduce the variables to only $N$ variables, namely
$Z=(z,w,u)\in\C^N$.
The main result of this section is the following.
\begin{Thm}\Label{gamma}
Under the assumptions of
Theorem~{\rm\ref{basic}}, there is an integer $s\ge 0$, a germ of a
holomorphic map
\begin{equation}\Label{}
\Gamma\colon \big(\C\times E^{2sl}\times\C^N,
\C\times E^{2sl}_0\times\{0\}\big)\to
(\C^N,0)
\end{equation}
with components in the ring
$\cR^{2sl}_N$,  and an integer
$r\ge 1$ such that for every $k$-equivalence $H$ between
$(M,0)$ and $(M',0)$ with $k> 2sl$, one has for
$Z=(z,w,u)$ sufficiently small,
\begin{equation}\Label{gam}
H(Z) =
\Gamma\Big(\frac{1}{q(j^{2sl}H(0,0,u))},
\D^{2sl}H(0,0,u),Z\Big)
+ O\Big(\frac{k-2sl}{r}\Big),
\end{equation}
where $q$ is the polynomial  $q^0_s$ on
$J^{2sl}_{N,N}$ given by Theorem~{\rm\ref{iterth}}.
\end{Thm}

\begin{Rem} {\rm The proof of Theorem \ref{gamma} shows that the integer
 $s\ge 0$ in this theorem
can be chosen to be the {\em Segre number} of $M$ at $0$
introduced in \cite{BER}. In particular, $s=0$ if and
only if $M$ is totally real, in which case the conclusion of
Theorem \ref{gamma} is obvious since $Z=u$.
In all other cases
we have $s\ge 1$.}
 \end{Rem}

Before proving Theorem \ref{gamma}, we shall state the
following corollary, which is of independent interest.
\begin{Cor}\Label{converge}
Under the assumptions of Theorem~{\rm\ref{basic}}
a formal equivalence $H$ between $(M,0)$ and $(M',0)$
is convergent if and only if the power series $j^{2sl}H(0,0,u)$ is
convergent in $u\in\C^{d_2}$.
\end{Cor}

\begin {proof}[Proof of Corollary {\rm\ref{converge}}]  Suppose that
$H$ is a formal equivalence.  If $H$ is convergent, it is clear
that $j^{2sl}H(0,0,u)$ is also convergent.  Conversely, if
$j^{2sl}H(0,0,u)$ is convergent, then the first term on the
right hand side of (\ref{gam}) is a convergent power series
in $Z$ by composition.  Since
$H$ is a $k$-equivalence for every $k$, the remainder term is
$0$, and hence $H(Z)$ is also convergent by (\ref{gam}).
\end{proof}

   For the proof of Theorem \ref{gamma}, we begin by defining
inductively a sequence of germs  of holomorphic maps
$$V^\k\colon
(\C^{\k n}\times\C^{d_2},0)\to(\C^N,0),
\quad \k=0,1,\ldots,$$
as follows. As before, we choose a holomorphic map $Q$
from a neighborhood of $0$ in
$\C^n\times\C^n\times\C^{d_1}\times\C^{d_2}$ to
$\bC^{d_1}$ satisfying (\ref{Qcond}) so that $M$ is given near $0$ by
(\ref{mq}). We put
$V^0(u):=(0,0,u)\in\C^n\times\C^{d_1}\times\C^{d_2}$ and
\begin{equation}\Label{Vind}
V^{\k+1}(t^0,t^1,\ldots,t^\k,u):=
\big(t^0,Q(t^0,\1{V^\k}(t^1,\ldots,t^\k,u)),u\big)
\in\C^n\times\C^{d_1}\times\C^{d_2}
\end{equation}
for $\k\ge 0$, $t^0,t^1,\ldots,t^\k\in\C^n$ and
$u\in\C^{d_2}$.
It is easy to check that for
$\k\ge 0$,
\begin{equation}\Label{vs1s}
\big(V^{\k+1}(t^0,t^1,\ldots,t^\k,u),
\1{V^\k}(t^1,\ldots,t^\k,u)\big)\in\3M,
\end{equation}
and hence also
\begin{equation}\Label{vss1}
\big(V^{\k}(t^1,\ldots,t^\k,u),
\1{V^{\k+1}}(t^0,t^1,\ldots,t^\k,u)\big)\in\3M.
\end{equation}
It will be convenient to introduce for every $s\ge 1$, the germ at $0$ of
a holomorphic map
\begin{multline}\Label{chis}\Xi^s(t^0,\ldots,t^{2s-1},u):=\\
\big(V^{2s}(t^0,\ldots,t^{2s-1},u),
\1{V^{2s-1}}(t^1,\ldots,t^{2s-1},u),
\ldots,\1{V^1}(t^{2s-1},u),V^0(u)\big).\end{multline}
Observe that the map
\begin{equation}\Label{}
\C^{2sn}\times\C^{d_2}\ni (t^0,\ldots,t^{2s-1},u)
\mapsto \Xi^s(t^0,\ldots,t^{2s-1},u)\in
\3M^{2s}\subset\C^{(2s+1)N}
\end{equation}
parametrizes a germ at $0$ of the submanifold of
$\3M^{2s}$ given by
$$\big\{(Z,\Z^1,Z^1,\ldots,\Z^s,Z^s)\in\3M^{2s} :
Z^s=(0,0,u)\big\}.$$
In this notation we have the following consequence of
Theorem~\ref{iterth}.
\begin{Cor}\Label{DDD}
Under the assumptions of
Theorem~{\rm\ref{basic}}, for any integer $s\ge 1$, there
exists a germ of a holomorphic map
\begin{equation}\Label{}
\Phi^s\colon \big(\C\times E^{2sl}\times\C^{2sn+d_2},
\C\times E^{2sl}_0\times\{0\}\big)\to
(\C^N,0)
\end{equation}
whose components are in the ring
$\cR^{2sl}_{2sn+d_2}$ such that,
if $H$ is a $k$-equivalence between $(M,0)$ and
$(M',0)$ with $k> 2sl$, then
\begin{multline}\Label{HV2s}
H(V^{2s}(t^0,\ldots,t^{2s-1},u))\equiv
\Phi^s\Big(\frac{1}{q^0_s(j^{2sl}H(0,0,u))},
\D^{2sl}H(0,0,u),t^0,\ldots,t^{2s-1},u\Big)
+ \\ r_H^{s}(t^0,\ldots,t^{2s-1},u),
\end{multline}
	where $q^0_s$ is the polynomial given by
Theorem~{\rm\ref{iterth}} and $r_H^s$ is a germ at $0$ of
a holomorphic map from $\C^{2ns+d_2}$ to $\C^N$
that vanishes of order $k-2sl$ at the origin.
\end{Cor}

\begin{proof}
We use (\ref{iter}) for $r=0$ and
substitute $\Xi^s(t^0,\ldots,t^{2s-1},u)$
for $(Z,\Z^1,Z^1,\ldots,\Z^s,Z^s)$, where $\Xi^s$ is given by
(\ref {chis}). The corollary easily follows by taking
$$\Phi^s(\T,\L,\L',z',t^0,\ldots,t^{2s-1},u):=
\Psi^{0,s}
\big(\T,\L,\L',z',\Xi^s(t^0,\ldots,t^{2s-1},u)\big)$$
and $r_H^s:=R^{0,s}\circ\Xi^s$.
\end{proof}

We next define a sequence of germs $v^\k$ at
$0$ of holomorphic maps from $\C^{\k n}$ to $\C^{n+d_1}$,
$\k\ge 0$, by
\begin{equation}\Label{vV}
V^s(t^0,\ldots,t^{\k-1},u)|_{u=0}=
\big(v^\k(t^0,\ldots,t^{\k-1}),0\big)\in
\C^{n+d_1}\times \C^{d_2}.
\end{equation}
Recall that the submanifold $M_0\subset\C^{n+d_1}$
defined by (\ref{mu}) is of finite type at $0$.
The map $v^\k$ defined above is the $\k$th {\em Segre
map} of $M_0$ in the sense of \cite{BER99}.
Hence by \cite{BER99} (Theorem 3.1.9)  the generic
rank of
$v^\k$ equals
$n+d_1$ for $\k$ sufficiently large. As in \cite{BER}
we call the smallest such $\k$ the {\em Segre number}
of $M_0$ at $0$ and denote it by s.
By \cite {BER99} (Lemma~4.1.3) we have
\begin{equation}\Label{dbl}
\max_{(x^1,\ldots,x^s)\in\3O}
\rk\frac{\d
v^{2s}}{\d(t^0,t^{s+1},t^{s+2},\ldots,t^{2s-1})}
(0,x^1,\ldots,x^{s-1},x^s,x^{s-1},\ldots,x^1)
= n+d_1.
\end{equation}
and
\begin{equation}\Label{equiv0}
v^{2s}(0,x^1,\ldots,x^{s-1},x^s,x^{s-1},\ldots,x^1)
\equiv 0,
\end{equation}
where $\3O$ is an arbitrary sufficiently small neighborhood of $0$ in $\C^{sn}$.
Note that in (\ref{dbl}) we differentiate only with
respect to the first vector $t^0$ and the last $s-1$
vectors $t^{s+1},\ldots,t^{2s-1}$.

For the proof of Theorem \ref{gamma}, we shall also
need the following special case of Proposition~4.1.18
in
\cite{BER99}.
\begin{Lem}
\Label{implicit}
Let
$$V\colon(\bC^{r_1}\times \C^{r_2},0)\to (\bC^N,0),
\quad r_2\ge N,
$$
be
a germ of a holomorphic map satisfying
$V(x,\xi)|_{\xi=0}\equiv 0$,  with
$(x,\xi)\in\C^{r_1}\times\C^{r_2}$, and 
for any sufficiently small neighborhood 
$\3O$ of $0$ in $\C^{r_1}$
\begin{equation}\Label{wid}
\max_{x\in \3O}\Big\{\rk\frac{\d
V}{\d\xi}(x,0) \Big\}=N.
\end{equation}
Then there exist germs of
holomorphic maps
\begin{equation}\Label{solution}
\delta\colon(\C^{r_1},0)\to \C, \quad
\delta(x)\not\equiv 0, \quad
\phi\colon(\C^{r_1}\times\C^N,0)\to(\C^{r_2},0)
\end{equation}
satisfying
\begin{equation}\Label{vz}
V\Big(x,\phi\big(x,\frac{Z}{\delta(x)}\big)\Big)
\equiv Z
\end{equation}
for all $(x,Z)\in \C^{r_1}\times\C^N$ such that $\delta(x)\ne
0$ and  both $x$ and $Z/\delta(x)$ are sufficiently small.
\end{Lem}

\begin{proof}[Proof of Theorem {\rm\ref{gamma}}]
We shall take $s$ to be the Segre number of $M_0$ at
$0$.  In the notation of Lemma
\ref{implicit} we take
$x=(x^1,\ldots,x^s)\in\C^{sn}$,
$\xi=(y,u)=(y^0,y^1,\ldots,y^{s-1},u)
\in\C^{sn}\times\C^{d_2}$
and set
\begin{equation}\Label{}
\begin{aligned}
L(x,y,u)&:=
(y^0,x^1,\ldots,x^s,x^{s-1}+y^{s-1},\ldots,x^1+y^1,u),
\\
V(x,\xi) &= V(x,y,u) := V^{2s}(L(x,y,u)),
\end{aligned}
\end{equation}
where $V^{2s}$ is defined by (\ref{Vind}). Here
$r_1:=sn$ and
$r_2:=sn+d_2$. Observe that $L$ is a linear automorphism of
$\C^{2sn+d_2}$. It follows from
(\ref{vV}) and (\ref{equiv0}) that $V(x,0)\equiv
0$. Furthermore it follows from (\ref{Vind}),
(\ref{vV}) and (\ref{dbl}) that condition
(\ref{wid}) also holds. Hence we can apply
Lemma~\ref{implicit}. Let
$$\delta\colon(\C^{sn},0)\to\C, \quad
\phi\colon (\C^{sn+N},0)\to(\C^{sn+d_2},0)$$
be given by the lemma, so that (\ref{vz}) holds.
By Corollary~\ref{DDD}, we obtain
\begin{multline}\Label{HV2s0}
H(Z)\equiv
\Phi^s
\bigg(
\frac{1}{q^0_s(j^{2sl}H(0,0,u))},
\D^{2sl}H(0,0,u),
L
\Big(x,
\phi\big(x,
\frac{Z}{\delta(x)}
\big)
\Big)
\bigg)
+\\
r_H^{s}
\bigg(L\Big(x,
\phi\big(
x,\frac{Z}{\delta(x)}
\big)\Big)\bigg),
\end{multline}
with $Z=(z,w,u)$.
By a simple change of $\Phi^s$ and $r^s_H$, we
obtain from (\ref{HV2s0}) the equivalent identity
\begin{equation}\Label{HV2s1}
H(Z)=
\2\Phi^s
\Big(
\frac{1}{q^0_s(j^{2sl}H(0,0,u))},
\D^{2sl}H(0,0,u),
\frac{Z}{\delta(x)},
x
\Big)
+
\2r_H^{s}
\Big(
\frac{Z}{\delta(x)},x
\Big)
\end{equation}
for all $(x,Z)\in\C^{ns+N}$ such that $\delta(x)\ne 0$
and both $x$ and $Z/\delta(x)$ are sufficiently small.
Here $\delta$ and  $\2\Phi^s$ are independent of $H$,
the  components of  $\2\Phi^s$ are in the ring
$\cR^{2sl}_{sn+N}$ and $\2r_H^{s}$ is a germ at $0$
in $\C^{sn+N}$, depending on $H$ and vanishing of
order $k-2sl$ at $0$.

Observe that the left-hand side of (\ref{HV2s1}) is
independent of the parameter $x\in\C^{ns}$, whereas
the right-hand side contains this parameter.
We choose $x_0\in\C^{ns}$ such that the function
$\2\l\mapsto\delta(\2\l x_0)$ does not vanish
identically for
$\2\l$ in a neighborhood of $0$ in $\C$, and put
$x=\2\l x_0$ in (\ref{HV2s1}). For convenience we
consider a holomorphic change of variable $\l=
h(\2\l)$ near the origin in $\C$, where $h$ is
determined by the identity
$\delta(\2\l x_0)=\l^m$ for an appropriate integer
$m\ge 0$. By a further simple change of $\2\Phi^s$ and
$\2r^s_H$, we conclude from (\ref{HV2s1}) that the identity
\begin{equation}\Label{HV2s2}
H(Z)\equiv
\widehat\Phi^s
\Big(
\frac{1}{q^0_s(j^{2sl}H(0,0,u))},
\D^{2sl}H(0,0,u),
\frac{Z}{\l^m},\l
\Big)+
\widehat r_H^{s}
\Big(
\frac{Z}{\l^m},\l
\Big),
\end{equation}
holds for all $(\l,Z)=(\l,z,w,u)\in\C^{1+N}$ such that $\l\ne 0$
and both $\l$ and $Z/\l^m$ are sufficiently small.
Again $\widehat\Phi^s$ is independent of $H$
and its components are in the ring
$\cR^{2sl}_{N+1}$ and $\widehat r_H^{s}$ is a germ at
$0$ in $\C^{N+1}$, depending on $H$ and vanishing of
order $k-2sl$ at $0$.

We next expand both sides of (\ref{HV2s2}) in Laurent series in $\l$
and equate the constant terms.
The required properties of those terms are established in the following lemma.
\begin{Lem}\Label{powers}
Let $V_0$ and $V_1$ be finite-dimensional vector
spaces with fixed linear coordinates
$x_0$ and $x_1$ respectively, and
$P(x_0,x_1,\l)$ be in the ring $\cR(V_{0}\times
V_{1}\times\C,V_0)$ with
$P(x_0,0,0)\equiv 0$. For  a fixed integer
$m\ge 0$, consider the Laurent series
expansion
$$P\big(x_0,\frac{x_1}{\l^m},\l\big)=
\sum_{\nu\in\bZ} c_\nu(x_0,x_1) \l^\nu. $$
Then $c_0(x_0,0)\equiv 0$ and, for every $\nu\in\bZ$,
$c_\nu$ is in the ring
$\cR(V_0\times V_1,V_0)$.
In addition, if $P=O(K)$ for some integer
$K>0$, then
$c_\nu=O(\frac{K-\nu}{m+1})$
for all $\nu\in\bZ$ such that $\nu\le K$.
\end{Lem}

\begin{proof}
We expand $P$ in power series of the form
\begin{equation}\Label{Pexp}
P(x_0,x_1,\l)=\sum P_{\b,\mu}(x_0)x_1^{\b}\l^\mu
=\sum P_{\a,\b,\mu}x_0^\a x_1^{\b}\l^\mu,
\end{equation}
with $P_{0,0}(x_0)\equiv 0$,
where $\a\in\bZ_+^{\dim V_0}$, $\b\in\bZ_+^{\dim
V_1}$, $\mu\in\bZ_+$.
Then $P_{\b,\mu}$ is a polynomial in $x_0$
satisfying the estimates (\ref{pa}).  Since
\begin{equation}\Label{sz}
c_\nu(x_0,x_1)=\sum_{\b,\mu\atop\mu - m|\b| = \nu}
P_{\b,\mu}(x_0)x_1^{\b}=
\sum_{\a,\b,\mu\atop\mu - m|\b| = \nu}
P_{\a,\b,\mu}x_0^\a x_1^{\b},
\end{equation}
we conclude that $c_0(x_0,0)\equiv P_{0,0}(x_0)\equiv
0$ and
$c_\nu\in\cR(V_0\times V_1,V_0)$ for every $\nu\in
\bZ$. Now assume that
$P=O(K)$. This means that
$\mu+|\a|+|\b|\ge K$ holds whenever
$P_{\a,\b,\mu}\ne 0$. For fixed $\nu\in\bZ$, this
inequality together with $\mu=\nu+m|\b|$ implies, in
particular, that
$\nu+(m+1)(|\a|+|\b|)\ge K$  in the last sum of
(\ref{sz}) or, equivalently,
$|\a|+|\b|\ge
\frac{K-\nu}{m+1}$ in that sum.
This completes the proof of the lemma.
\end{proof}

We now complete the proof of Theorem~\ref{gamma}.
We expand the right-hand side of (\ref{HV2s2})
in Laurent series in $\l$. Since the left-hand side
is independent of $\l$, we equate it to the constant
term of the Laurent series.
The required conclusion of Theorem~\ref{gamma} with $r:=m+1$ follows
by applying Lemma~\ref{powers} for $\nu=0$ to
$\widehat\Phi^s$ with $V_0:=\C\times E^{2sl}_0$,
$V_1=E^{2sl}_1\times \C^N$ and to
$\widehat r^s_H$ with $V_0:=0$, $V_1:=\C^N$.  The
proof of Theorem \ref{gamma} is now complete.
\end{proof}

\section{Equations in jet spaces}\Label{univsec}
In Theorem~\ref{gamma} we showed that
every $k$-equivalence $H$ between $(M,0)$ and $(M',0)$
for $k$ sufficiently large satisfies the identity (\ref{gam}),
i.e. is parametrized up to the given order by the jet $j^{2sl}H(0,0,u)$.
However, it will be more convenient to regard
$\Theta_H(u):=\big(\frac{1}{q(j^{2sl}H(0,0,u))},\D^{2sl}H(0,0,u)\big)$
as the main parameter since the parametrization then
becomes polynomial rather than rational.
Our goal in this section is to give a set of
equations such that any germ at $0$ of a
real-analytic map
$\bR^{d_2}\ni u\mapsto\Theta(u)\in \C\times E^{2sl}$
satisfies these equations if and only if the mapping
$\bC^N\ni(z,w,u)\mapsto
\Gamma(\Theta(u),z,w,u)\in\bC^N$ is a germ at $0$ of
a holomorphic self map of $\bC^N$ sending
$M$ into
$M'$; here $\Gamma$ is the mapping given by
(\ref {gamma}).

We shall need a real analogue of the ring
$\cR(V,V_0)$ defined in \S\ref{rings}.
Given a finite dimensional {\em real}
vector space $W$  and a real vector subspace
$W_0\subset W$, define $\cR_\bR(W,W_0)$ to be the
ring of all germs of real valued real-analytic
functions $f$ at
$W_0$ in $W$ such that  all partial derivatives
$\d^\alpha f|_{W_0}$ are real polynomial functions
on $W_0$. In the following we shall consider the
spaces $\C$ and $E^{2sl}$ as real vector spaces
and real-analytic functions on these spaces
with respect to real and imaginary parts of
vectors in these spaces.

\begin{Thm}\Label{universal}
Assume that the conditions of Theorem~{\rm\ref{gamma}}
are satisfied, and let
$s$, $q$, and
$\Gamma$ be given by  that
theorem. Then there exist
a finite collection of functions
$f_j\in\cR_\bR\big(\C\times E^{2sl}\times\bR^{d_2},
\C\times E^{2sl}_0\times\{0\}\big)$,
$1\le j\le j_0$,
and positive real numbers $a$ and $b$, with $b\ge
2sla$, such that the following hold.
\begin{enumerate}
\item[(i)]
For every $k$-equivalence $H$ between
$(M,0)$ and $(M',0)$ with $k> b/a$, one has
\begin{equation}\Label{fs}
f_j\Big(\frac{1}{q(j^{2sl}H(0,0,u))},
\D^{2sl}H(0,0,u),u\Big)=O(|u|^{ak-b}),\quad
1\le j\le j_0.
\end{equation}
\item[(ii)] For every germ
$\Th\colon(\bR^{d_2},0)\to (\C\times
E^{2sl},\C\times E^{2sl}_0)$
of a real-analytic map satisfying
\begin{equation}\Label{approx}
f_j(\Theta(u),u)\equiv 0,\quad
1\le j\le j_0,
\end{equation}
the germ
$\Gamma_\Th\colon (\C^n\times\C^{d_1}\times\bR^{d_2},0)\to (\C^N,0)$
of the real-analytic map defined by
\begin{equation}\Label{Hth}
\Gamma_\Th(z,w,u):= \Gamma(\Th(u),z,w,u),
\end{equation}
extends to a germ at $0$ of a holomorphic map of
$\C^N$ into itself sending
$(M,0)$ into $(M',0)$.
\end{enumerate}
\end{Thm}

\begin{Rem}
{\rm It should be mentioned that the holomorphic
extension of the germ
$\Gamma_\Th$ defined by {\rm(\ref{Hth})} need not be
invertible.}
\end{Rem}

Before starting the proof of Theorem~\ref{universal}
we shall need a composition lemma for the rings
$\cR_\bR(W,W_0)$ whose complex analogue is
a special case of Lemma~\ref{compose}.

\begin{Lem}\Label{rcompose}
Let $W_0$, $W_1$ and $\2W$ be finite-dimensional
real vector spaces with fixed bases. Denote by
$x_0$, $x_1$ and $\2x$ the corresponding real linear
coordinates in these spaces. Let
$\phi\colon(W_0\times W_1,W_0)\to(\2W,0)$
be a germ at $W_0$ of a real-analytic map
whose components are in the ring
        $\cR_\bR(W_0\times W_1,W_0)$.
Then, for every germ
$\2f\colon (\2W,0)\to \bR$ of a real-analytic map,
there exists $f\in\cR_\bR(W_0\times W_1,W_0)$  such
that
$f(x_0,x_1)=\2f(\phi(x_0,x_1))$.
\end{Lem}

Lemma~\ref{rcompose} is a straightforward consequence
of the chain rule and is left to you, gentle reader.

\begin{proof}[Proof of Theorem~{\rm\ref{universal}}]
We continue to work with normal coordinates near the
origin $Z=(z,w,u)$ for $M$.
We fix a local parametrization
of $M$ at $0$ of the form
$$\bR^{2n+d_1}\times\bR^{d_2}\ni (t,u)\mapsto
(z(t,u),w(t,u),u)\in M\subset\C^N.$$
Let
$\rho'(Z',\1{Z'})=\big(\rho^1{}'(Z',\1{Z'}),\ldots,\rho^d{}'(Z',\1{Z'})\big)$
be a defining
function for $M'$ near $0$
and $\Gamma$ be given by Theorem~\ref{gamma}.
For $1\le i\le d$, $\a\in\bZ_+^{2n+d_1}$
and $\Th\in \C\times E^{2sl}$,
we consider the functions
\begin{equation}\Label{fja}
f^i_{\a}(\Th,u):=\frac{\d}{\d t^\a}
\rho^i{}'
\Big(\Gamma\big(\Th,z(t,u),w(t,u),u\big),
\ov{\Gamma\big(\Th,z(t,u),w(t,u),u\big)}
\Big)\Big|_{t=0}.
\end{equation}
It follows from the properties of $\Gamma$,
the chain rule and Lemma~\ref{rcompose}
that the $f^i_{\a}$ are
in the ring $\cR_\bR\big(\C\times
E^{2sl}\times\bR^{d_2},
\C\times E^{2sl}_0\times\{0\}\big)$.
If follows from the definition that we can think of
this ring as a subring of the following formal power
series ring with polynomial coefficients
\begin{equation}\Label{fring}
\bR\big[\Re\T,\Im\T,\Re\L,\Im\L\big]
\big[\!\big[\Re\L',\Im\L',\Re z',\Im
z',u\big]\!\big],
\end{equation}
where $(\T,\L,\L',z')$ are complex coordinates in
$$\C\times J^{2sl}_{N,N}\times J^{2sl}_{n,d}\times
\C^n=\C\times E^{2sl}.$$
It is a standard fact from commutative algebra
that any formal power series ring with coefficients
in a Noetherian ring is again Noetherian;
in particular, the ring (\ref{fring}) is Noetherian.
Hence there exists an integer $m_0\ge 0$ such that
the subset
\begin{equation}\Label{set}
\big\{f^i_{\a}:1\le i\le d, |\a|\le m_0\big\}
\end{equation}
generates the same ideal in the ring (\ref{fring})
as all the $f^i_{\a}$, $1\le i\le d$,
$\a\in\bZ_+^{2n+d_1}$.

By the identity (\ref{gam}) and the definition of
$k$-equivalence we have for $1\le i\le d$ and $|\a|\le
m_0$,
\begin{multline}\Label{}
f^i_{\a}\Big(\frac{1}{q(j^{2sl}H(0,0,u))},
\D^{2sl}H(0,0,u),u\Big)= \\
\frac{\d}{\d t^\a}
\rho^i{}'
\Big(H\big(z(t,u),w(t,u),u\big), \ov{H\big(z(t,u),w(t,u),u\big)}
\Big)\Big|_{t=0}
+ O\Big(\frac{k-2sl}{r}-m_0\Big)
= \\ O(k-m_0) + O\Big(\frac{k-2sl}{r}-m_0\Big).
\end{multline}
Hence we proved (\ref{fs})
with the collection $f_j$, $1\le j\le j_0$,
being the set of functions given by (\ref{set})
and $a:=1/r\le 1$, $b:=m_0+(2sl/r)$.
This completes the proof of (i).

We shall now prove  (ii).
By the choice of the set (\ref{set}),
every germ $f^i_{\a}(\Th,u)$ can be written in the
form
\begin{equation}\Label{genid}
f^i_{\a}(\Th,u)\equiv
\sum_{j=1}^{j_0}c_j(\Th,u)f_j(\Th,u),
\end{equation}
where $c_j(\Th,u)$  are in the ring given by
(\ref{fring}). Since $\Th(0)\in\C\times E^{2sl}_0$,
the germ $\Th(u)$ can be substituted for $\Th$ in
each $c_j(\Th,u)$  to obtain a formal
power series in
$\bR[[u]]$. From  (\ref{genid}) and the assumption
(\ref{approx}) on $\Th(u)$ we obtain the following
identities of convergent power series in $u$:
$$f^i_{\a}\big(\Th(u),u\big)\equiv 0,
\quad 1\le i\le d, \quad \a\in\bZ_+^{2n+d_1}.$$
In view of (\ref{fja}) we conclude
that
$\rho'(\Gamma_\Th(z,w,u),\1{\Gamma_\Th(z,w,u)})= 0$
for $(z,w,u)\in M$ near the origin. This
completes the proof of (ii) and hence that of
Theorem~\ref{universal}.
\end{proof}

\section{Artin and Wavrik theorems}

We state two approximation results due to Artin
\cite{artin1} and Wavrik \cite{Wav} which will be used
(in conjunction with Theorem~\ref{universal}) in the
proof of Theorem~\ref{tech}. We start by stating the
result of Artin, which  implies that any formal
solution of a system of real-analytic equations may
be approximated to any preassigned order by a
convergent solution of that system. We use the
superscripts $f$, $c$, $a$ to denote formal,
convergent and approximate solutions respectively.

\begin{Thm}[\cite{artin1}]\Label{artin}
Let $g_j(t,u)\in\bR\{t,u\}$, $1\le j\le j_0$, be
convergent power series in $t=(t_1,\ldots,t_\delta)$
and $u=(u_1,\ldots,u_\gamma)$. Then for any integer
$\k\ge1$ and any formal power series
$t^f(u)\in\big(\bR[[u]]\big)^\delta$, satisfying
$$t^f(0)=0, \quad g_j(t^f(u),u)\equiv 0, \quad 1\le
j\le j_0,$$ there exists a convergent power series
$t^c(u)\in\big(\bR\{u\}\big)^\delta$ satisfying
$$t^c(u)=t^f(u)+O(|u|^\k), \quad
g_j(t^c(u),u)\equiv 0, \quad 1\le
j\le j_0.$$
\end{Thm}

We now turn to a result of Wavrik
which states that an approximate formal solution
of a system of formal equations of
a certain type may be approximated by
an exact formal solution of that system.
(This result generalizes another result of
Artin~\cite{artin2} which deals with more special
systems of equations. See also
Denef-Lipshitz \cite{DL} for related results.)

\begin{Thm}[\cite{Wav}]\Label{wavrik}
Let
$h_j(x,y,u)\in\bR\big[x\big]\big[\!\big[y,u\big]\!\big]$,
$1\le j\le j_0$, be formal power series in
$y=(y_1,\ldots,y_\beta)$ and
$u=(u_1,\ldots,u_\gamma)$  with coefficients which are
polynomials in $x=(x_1,\ldots,x_\a)$. Then for any
integer $\k\ge 1$, there exists an integer
$\eta\ge 1$ such that, for any formal power series
$x^a(u)\in\big(\bR[[u]]\big)^\a$,
$y^a(u)\in\big(\bR[[u]]\big)^\b$ satisfying
$$y^a(0)=0, \quad h_j(x^a(u),y^a(u),u)=O(|u|^\eta),
\quad 1\le j\le j_0,$$
there exist formal power series
$x^f(u)\in\big(\bR[[u]]\big)^\a$,
$y^f(u)\in\big(\bR[[u]]\big)^\b$ satisfying
\begin{equation}\Label{xyh}
\begin{aligned}
  x^f(u)=x^a(u)+O(|u|^\k), \quad
y^f(u)=y^a(u)+O(|u|^\k), \\
  h_j(x^f(u),y^f(u),u)\equiv 0, \quad 1\le
j\le j_0.\quad\quad \quad \quad \quad
\end{aligned}
\end{equation}
\end{Thm}
An immediate corollary of Theorems \ref{artin} and
\ref{wavrik}, which we shall need, is the following.
\begin{Cor}\Label{awcor}
Let $X$, $Y$, $U$ be real finite-dimensional vector
spaces with fixed linear coordinates
$x=(x_1,\ldots,x_\a)$, $y=(y_1,\ldots,y_\b)$,
$u=(u_1,\ldots,u_\gamma)$ respectively
and let $h_j(x,y,u)$, $1\le j\le
j_0$, be germs of functions in the ring
$\cR_\bR(X\times Y\times U,X)$.
Then for
any integer $\k\ge 1$, there exists an integer
$\eta\ge 1$ such that, for any germs at $0$ of
real-analytic maps $x^a\colon (U,0)\to X$,
$y^a\colon (U,0)\to Y$ satisfying
$$y^a(0)=0, \quad h_j(x^a(u),y^a(u),u)=O(|u|^\eta),
\quad 1\le j\le j_0,$$
there exists germs at $0$ of real-analytic maps
$x^c\colon (U,0)\to X$,
$y^c\colon (U,0)\to Y$ satisfying
$$x^c(u)=x^a(u)+O(|u|^\k), \quad
y^c(u)=y^a(u)+O(|u|^\k), \quad
h_j(x^c(u),y^c(u),u)\equiv 0, \quad 1\le
j\le j_0.$$
\end{Cor}

\begin{proof}
Let $h_j\in \cR_\bR(X\times Y\times U,X)$, $1\le
j\le j_0$, be given. It follows from the definition
of the ring $\cR_\bR(X\times Y\times U,X)$
that $h_j$ can be viewed as an element of
$\bR\big[x\big]\big[\!\big[y,u\big]\!\big]$. By
Theorem~\ref{wavrik}, given $\k\ge 1$, there exists
$\eta\ge 1$ such that if
$x^a$ and $y^a$ are as in the corollary,
there exist $x^f(u)\in\big(\bR[[u]]\big)^\a$,
$y^f(u)\in\big(\bR[[u]]\big)^\b$ satisfying
(\ref{xyh}). We may now apply Theorem~\ref{artin} with
$t=(x-x^a(0),y)$ (and hence $\delta=\a+\b$)
to conclude that there exists
$t^c(u)=(x^c(u)-x^a(0),y^c(u))$
with $x^c$, $y^c$ satisfying the conclusion of the
corollary.
\end{proof}

\section{End of proof of
Theorem~\ref{tech}}\Label{final} We keep the
notation used in \S\ref{univsec} and, in
particular, that of Theorem~\ref{universal}. Let
$s$, $q$, and
$\Gamma$  be given by Theorem~\ref{gamma}, and let
$H$ be a
$k$-equivalence between
$(M,0)$ and
$(M',0)$  with $k> 2sl$.
Define the germ
$\Theta_H\colon (\bR^{d_2},0)\to(\C\times
E^{2sl},\C\times E^{2sl}_0)$ of a real-analytic map by
$$\Th_H(u):= \left(\frac{1}{q(j^{2sl}H(0,0,u))},
\D^{2sl}H(0,0,u)\right).$$
Recall that $q(j^{2sl}H(0,0,0))\ne 0$
by (\ref{detprop}).
Let $f_j$, $1\le j\le j_0$, be given by
Theorem~\ref{universal}. By
part (i) of that theorem we have for $k> b/a\ge 2sl$,
\begin{equation}\Label{ak-b}
f_j(\Th_H(u),u)=O(|u|^{ak-b}), \quad 1\le j\le j_0.
\end{equation}
We shall now apply Corollary~\ref{awcor} to the system
of equations
$f_j(\Th (u),u)=0$, $1\le j\le j_0$.
Indeed, if $\Th\colon(\bR^{d_2},0)\to \big(\C\times
E^{2sl},\C\times E^{2sl}_0\big)$
is a germ of a real-analytic map, we may write
$\Th(u)=
(x(u),y(u))\in (\C\times E^{2sl}_0)\times
E^{2sl}_1$ so that the system of equations above
becomes $f_j(x(u),y(u),u) = 0$, $1\le j\le j_0$. Given
$\k>1$ we conclude by Corollary~\ref{awcor} that
there exists
$\eta\ge 1$ such that, if $ak-b\ge\eta$, the identity
(\ref{ak-b}) implies the existence of a germ
$\Th^c\colon(\bR^{d_2},0)\to\big(\bC\times
E^{2sl},\C\times E^{2sl}_0\big)$ of a  real-analytic
map satisfying
\begin{equation}\Label{lastid}
\Th^c(u)=\Th_H(u)+O(|u|^\k), \quad
f_j(\Th^c(u),u)\equiv 0, \quad
1\le j\le j_0.
\end{equation}
Then, by Theorem~\ref{universal}~(ii),
we conclude that
$\Gamma_{\Th^c}$ defined by (\ref{Hth}) extends as a
germ of a holomorphic map
$ (\C^N,0)\to (\C^N,0)$
  sending
$(M,0)$ into $(M',0)$.  We take
$\4H(Z):=\Gamma_{\Th^c}(Z)$. The first
identity in (\ref{lastid})
implies
\begin{equation}\Label{overlast}
\4H(Z)=\Gamma_{\Th_H}(Z)+O(|Z|^\k).
\end{equation}
On the other hand, by Theorem \ref{gamma} and,
in particular,  (\ref{gam}), we have for $k >2sl$
\begin {equation}\Label{last}
  H(Z)-\Gamma_{\Th_H}(Z)=
O\left({k-2sl\over r}\right).
\end{equation}
By increasing
$k$ if necessary, we can assume that $(k-2sl)/r\ge
\k$   so that $\4H(Z) = H(Z) + O(|Z|^\k)$ by
(\ref{overlast}) and (\ref{last}). Since $H$ is
invertible and
$\k > 1$, it follows that $\4H$ is also invertible.
  The proof of Theorem~\ref{tech} is now complete.

\section{CR equivalences}\Label{rfinal}  

If $M$ and $M'$ are real-analytic CR submanifolds of
$\bC^N$, with $p\in M$ and $p'\in M'$,  and $h: (M,p)\to
(M',p')$ is a germ of a mapping of class
$C^k$,
$1\le k\le
\infty$, we say that
$h$ is a {\em germ of a CR map of class $C^k$} if the
differential of
$h$ sends the CR bundle of $M$ into  that of $M'$.  If, in
addition, $h$ is a diffeomorphism at $p$ we shall say that $h$
is a {\em CR equivalence of class $C^k$} between $(M,p)$ and
$(M',p')$.
It is standard that the Taylor power series of any CR
equivalence of class $C^\infty$ between
$(M,p)$ and $(M',p')$ induces a formal
equivalence between $(M,p)$ and $(M',p')$.  Similarly, the
$k$th Taylor polynomial of any  CR equivalence of class $C^k$
induces a $k$-equivalence (see e.g. \cite{BER},
Proposition~1.7.14). Hence 
Corollary~\ref{main} implies the following.

\begin{Cor}\Label{nmain}
Let $M\subset \bC^N$ be a connected
real-analytic CR submanifold. Then there exists a
closed, proper real-analytic subvariety $V\subset M$ such
that for every
$p\in M\setminus V$, every real-analytic CR submanifold
$M'\subset
\bC^N$, and every $p'\in
M'$, the following are equivalent:
\begin{enumerate}
\item [(i)] $(M,p)$ and $(M',p')$ are $k$-equivalent
for all
$k>1$;
\item [(ii)] $(M,p)$ and $(M',p')$ are CR equivalent
of class $C^k$ for all finite
$k>1$;
\item [(iii)] $(M,p)$ and $(M',p')$ are formally
equivalent;
\item [(iv)] $(M,p)$ and $(M',p')$ are CR
equivalent of class $C^\infty$;
\item [(iii)] $(M,p)$ and $(M',p')$ are
biholomorphically equivalent.
\end{enumerate}
\end{Cor}

\end{document}